\documentstyle[11pt]{article}

\newcommand{\qed}{\hfill$\Box$}

\newcommand{\inv}{{\it inv}}
\newcommand{\Des}{{\it Des}}
\newcommand{\D}{{\it Des}}
\newcommand{\des}{{\it des}}

\newcommand{\Del}{{\it Del}}
\newcommand{\del}{{\it del}}
\newcommand{\maj}{{\it maj}}
\newcommand{\rmaj}{{\it rmaj}}

\newcommand{\supp}{{\it supp}}

\newtheorem{thm}{Theorem}[section]
\newtheorem{pro}[thm]{Proposition}
\newtheorem{lem}[thm]{Lemma}

\newtheorem{cor}[thm]{Corollary}
\newtheorem{fac}[thm]{Fact}
\newtheorem{obs}[thm]{Observation}

\newtheorem{df}[thm]{Definition}
\newtheorem{rem}[thm]{Remark}
\newtheorem{nt}[thm]{Note}

\def\ep{\epsilon}
\def\sg{\sigma}

\begin{document}

\title{Permutation Statistics on the Alternating Group }

\author{Amitai Regev\thanks{Partially supported by Minerva
Grant No. 8441 and by EC's IHRP Programme, within the Research Training Network ``Algebraic Combinatorics
in Europe'', grant HPRN-CT-2001-00272}\\
Department of Mathematics\\
The Weizmann Institute of Science\\
Rehovot 76100, Israel\\
{\em regev@wisdom.weizmann.ac.il}\\
\and
Yuval Roichman\thanks{Partially supported by the Israel Science Foundation, founded by the Israel Academy of Sciences and Humanities and by EC's IHRP Programme, within the Research Training Network ``Algebraic Combinatorics
in Europe'', grant HPRN-CT-2001-00272}\\
Department of Mathematics\\
Bar Ilan University\\
Ramat Gan 52900, Israel\\
{\em yuvalr@math.biu.ac.il}}

\maketitle

\bibliographystyle{acm}

\begin{abstract}
Let $A_n\subseteq S_n$ denote the alternating and 
the symmetric groups on $1,\ldots ,n$. MacMahaon's theorem~\cite{MM}, 
about the equi-distribution of the length and the major indices in $S_n$, 
has received far reaching refinements and generalizations, by 
Foata~\cite{F}, Carlitz~\cite{Ca1, Ca2},
Foata-Sch\"utzenberger~\cite{FS}, Garsia-Gessel~\cite{GG}
and followers.
Our main goal is to find analogous statistics and 
identities for the alternating group $A_{n}$.
A new statistic for $S_n$, {\it the delent number}, is introduced.
This new statistic is involved with new $S_n$ equi-distribution
 identities, refining 
some of the results in~\cite{FS} and~\cite{GG}. By a certain covering
map $f:A_{n+1}\to S_n$, such $S_n$ identities are `lifted' to $A_{n+1}$,
yielding the corresponding $A_{n+1}$ equi-distribution identities.
\end{abstract}

\section{Introduction}

\subsection{General outline}

One of the most active branches in 
enumerative combinatorics is the study of {\em permutation statistics}.
Let $S_n$ be the symmetric group on $1,\dots,n$. One is interested in the  refined count of permutations according to (non-negative, integer valued) combinatorial parameters.
For example, the number of inversions in a permutation - 
namely its {\it length} - is such a parameter. 
Another parameter is MacMahon's {\it major index}, 
which is defined via the {\it descent} set of a permutation - see below.

Two parameters that have the same generating function are said to be 
{\em equi-distributed}.
Indeed, MacMahon~\cite{MM} proved the remarkable fact
that the inversions and the major-index
statistics are equi-distributed on $S_n$.
MacMahon's classical theorem~\cite{MM}  has received far reaching
refinements and generalizations, including: 
multivariate refinements which imply equi-distribution on certain subsets of 
permutations 
(done by Carlitz~\cite{Ca1, Ca2}, Foata-Sch\"utzenberger~\cite{FS} and
Garsia-Gessel~\cite{GG});
analogues for other combinatorial objects, cf.~\cite{F, K, St};
generalizations to other classical Weyl groups, cf.~\cite{Rei, AR, ABR}.

Let $A_n\subseteq S_n$ denote the alternating group
on $1,\ldots ,n$.
Easy examples show that the above statistics fail to be
equi-distributed when restricted to $A_n$.
Our main goal is to find statistics on $A_n$ which are natural 
generalizations of the $S_n$ statistics and are equi-distributed on $A_n$, 
yielding analogous identities for their generating functions.
This goal is achieved by proving further refinements of the above
 $S_n$-identities.

It is well known that the above
statistics on $S_n$ may be defined via the Coxeter generators
$\{ (i,i+1)\mid 1\le i \le n-1\}$ of $S_n$.
Mitsuhashi~\cite{Mi}  pointed out at a certain set of generators
of the alternating group $A_n$, which play a role similar
to that of the above  Coxeter generators of $S_n$, see Subsection 1.3 below.
We use these generators to define the analogous length and descent 
statistics on the alternating group.

The $S_n$-Coxeter generators allow one to introduce the classical
canonical presentation of the elements of $S_n$, see Subsection 3.1. 
Similarly, the above Mitsuhashi's `Coxeter' generators
allow us to introduce the corresponding canonical presentation of the 
elements of $A_{n+1}$, see Subsection 3.3.
We remark that usually, $S_n$ is viewed as a double cover of $A_n$. However,
the above canonical presentations enable us to introduce
a covering map $f$ from the alternating group $A_{n+1}$ onto $S_n$, and thus
$A_{n+1}$ can be viewed as a covering of $S_n$.  

A new statistic, the {\it delent number}, plays a crucial role
in the paper, and allows us to `lift' $S_n$ identities to $A_{n+1}$.
The delent number on $S_n$ 
may be defined as follows: if the transposition
$(1,2)$ appears $r$ times in the canonical presentation of 
$\sg\in S_n$ then the delent number of $\sg$, $del_S(\sg)$, is $r$. 
An analogous statistic is defined for $A_{n+1}$, see Definition~\ref{dfDEL}.
We give direct combinatorial characterizations of this statistic
(see Propositions~\ref{pro13-i} and~\ref{pro12-i}  below) and show that
this statistic is involved in 
new $S_n$ equi-distribution identities,  
refining some of the results of Foata-Sch\"utzenberger~\cite{FS} 
and of Garsia-Gessel~\cite{GG}.
Identities involving the delent number
are then `lifted' by the covering map $f$, yielding
$A_{n+1}$ equi-distribution identities, see Theorem~\ref{thm101}, Theorem~\ref{shuff.1} 
and Corollary~\ref{shuff.2}.

In the Appendix we present different statistics on $A_n$,
and a consequent different analogue of MacMahon's equi-distribution
theorem. These statistics
are compatible with the usual point of view of $S_n$ 
as a double cover of $A_n$.

The above setting and results are connected with enumeration 
of other combinatorial objects, such as permutations avoiding patterns,
leading to $q$-analogues of the classical $S_n$ statistics
and of the Bell and Stirling numbers. A detailed study of these $q$-analogues is given 
in~\cite{RR} (a few of these results appear in 
Subsection 5.3).

\bigskip

The paper is organized as follows :
The rest of this section surveys briefly the classical background
and lists our main results. 
 Background and notations are 
given in detail in Section 2, while the $A$ canonical presentation is analyzed in
Section 3. In Section 4 we study the length statistics, and in Section 5 
we discuss the relations between various $S$- and $A$-statistics, 
relations given by the
map $f:A_{n+1}\to S_n$. In Section 6 we study the ordinary and the 
reverse major indices, together with the delent statistics. Additional 
properties of the delent numbers are given in Section 7. In Section 8 
we prove some lemmas on shuffles - lemmas that are needed for the 
proof of the main theorem.
The main theorem (Theorem~\ref{shuff.1}) and its proof are given in 
Section 9. Finally, the Appendix constitutes Section 10.

\subsection{Classical $S_n$-Statistics}\label{class}

Recall that the Coxeter generators  $S:=\{(i,i+1)|\ 1\le i\le n-1\}$
of $S_n$ give rise to various combinatorial statistics, like the {\it length}
statistic, etc. As we show later, most of these $S_n$ statistics have
$A_n$ analogues, therefore we add ``S-'' and ``A-'' to the titles of 
the corresponding statistics.\\  
{\it The S-length}: For $\pi\in S_n$
let $\ell _S(\pi)$ be the standard length of $\pi$ with respect to
these Coxeter generators.\\
{\it The S-descent}: 
Given a permutation $\pi$ in the symmetric group 
$S_n$, the {\it $S$-descent set} of $\pi$ is defined by
$$
\D _S(\pi):=\{i\ |\ \ell_S(\pi)>\ell_S(\pi s_i)\}
=\{i\mid \pi (i) > \pi (i+1)\}.
$$

\noindent
The {\it descent number} of $\pi$, $\des _S(\pi)$, is defined by
$\des _S(\pi):=|\D _S(\pi)|.$

\noindent
The {\it major index}, $\maj _S(\pi)$ is 
$$ 
\maj _S(\pi) := \sum\limits_{i\in \D _S(\pi)} i.
$$
The corresponding  {\it reverse major index} does depend on $n$, and is denoted
$$
\rmaj _{S_n}(\pi):=\sum\limits_{i\in\D _S(\pi)}(n-i).
$$
The reverse major index $rmaj_{S_n}(\pi)$ is implicit in~\cite{FS}.

\medskip

These statistics are involved in many combinatorial identities.
First, MacMahon proved the following equi-distribution of the length and 
the major indices \cite{MM}:
\[
\sum_{\sg\in S_n}q^{\ell _S (\sg)}=
\sum_{\sg\in S_n}q^{maj _S (\sg)}.
\]
Foata~\cite{F} gave a bijective proof of MacMahon's theorem, then 
Foata and Sch\"utzenberger~\cite{FS} applied this bijection  to refine MacMahon's identity
by analyzing bivariate distributions.
Garsia and Gessel~\cite{GG} extended the analysis to multivariate distributions.
Extensions of MacMahon's identity to hyperoctahedral groups appear in \cite{ABR}.

Combining Theorems 1 and 2 of~\cite{FS} 
one deduces the identity

\begin{thm}\label{class1}
For any subset $D_1\subseteq \{1,\ldots ,n-1\}$
$$
\sum\limits_{\{\pi\in S_n|\ \D_S(\pi^{-1})\subseteq D_1\}} q^{\maj_{S_n}(\pi)}
=
\sum\limits_{\{\pi\in S_n|\ \D_S(\pi^{-1})\subseteq D_1\}} q^{\rmaj_{S_n}(\pi)}
$$
$$
=\sum\limits_{\{\pi\in S_n|\ \D_S(\pi^{-1})\subseteq D_1\}} q^{\ell_S(\pi)}.
$$
\end{thm}

A bivariate equi-distribution follows.

\begin{cor}\label{class2} 
$$
\sum\limits_{\pi\in S_n} q_1^{\maj_{S_n}(\pi)}
q_2^{\des_S(\pi^{-1})}=
\sum\limits_{\pi\in S_n} q_1^{\rmaj_{S_n}(\pi)}
q_2^{\des_S(\pi^{-1})}=
\sum\limits_{\pi\in S_n} q_1^{\ell_S(\pi)}
q_2^{\des_S(\pi^{-1})}.
$$
\end{cor}

As already mentioned, one of the main goals in this paper is to find 
analogous statistics and 
identities for the alternating group $A_n$.
In the process we first prove some further refinements of some of the above
identities for $S_n$, refinements involving the new {\it delent} statistic,
see Theorems~\ref{thm101}.1  and~\ref{shuff.1}.1.

\subsection{Main Results}

Here is a summary of the main results of this paper.

\subsubsection{$A_n$-Statistics}\label{stat}

Following Mitsuhashi~\cite{Mi} we let 
$$
a_i:=s_1s_{i+1}=(1,2)(i+1,i+2)\qquad (1\le i\le n-1).
$$
Thus $a_i=a_i^{-1}$ if $i\ne 1$, while $a_1^2=a_1^{-1}$.
The set $A:=\{a_i\ | \ 1\le i\le n-1\}$
generates the alternating group on $n+1$ letters $A_{n+1}$ 
(see e.g. \cite{Mi}).
It is the above exceptional property of $a_1$ among the elements of $A$ - 
which naturally leads to the `delent' statistic (Definition 1.5 below), both 
for $S_n$ and for $A_{n+1}$. This new statistic enables us to deduce 
new refinements of the MacMahon-type identities for $S_n$, and for 
each such an identity to derive the analogous identity 
for $A_{n+1}$.


\medskip

The canonical presentation in $S_n$ by the Coxeter generators is well known,
and is discussed in Section~\ref{capre}, see Theorem~\ref{thm1}. 
With the above generating set $A$ of $A_{n+1}$ we also have 
canonical presentations for the elements 
of   $A_{n+1}$, as follows.
For each $1\le j\le n-1$ define 
\begin{eqnarray}\label{eqn2}
R^A_j=\{ 1,\;a_j,\;a_ja_{j-1},\;\ldots ,\; a_j\cdots a_2,\;a_j
\cdots a_2a_1,\;a_j\cdots a_2a^{-1}_1\}.
\end{eqnarray}

\begin{thm}\label{thm2i} (See Theorem \ref{thm2})
Let $v\in A_{n+1}$, then there exist unique elements $v_j\in R^A_j$,
 $1\le j\le n-1$, such that $v=v_1\cdots v_{n-1}$, and this presentation 
is unique.
Call that presentation $v=v_1\cdots v_{n-1}$ the {\it $A$ canonical 
presentation} of $v$. 
\end{thm}

\medskip

The $A$ canonical presentation 
allow us to introduce the A-length of an element in $A_{n+1}$:
\begin{df}\label{df1}
Let $v\in A_{n+1}$ with $v=a^{\ep _1}_{i_1}\cdots a^{\ep _r} _{i_r}$
($\ep_i=\pm 1$) its $A$ canonical presentation,
then its A-length is $\ell_A(v)=r$.

\end{df}

%
A combinatorial interpretation of the $A$-length in terms of 
inversions is given below, see Proposition~\ref{pro13}.

 The $A$-descent statistic is defined using the above generating set $A$ :

\begin{df}\label{df8} 
\begin{enumerate}
\item
The {\it alternating-descent} (i.e.~the $A$-descent) set
 of $\sigma\in A_{n+1}$ is defined by:
$$
\D_A(\sigma):=\{ 1\le i\le n-1\mid \ell_A(\sigma)\ge \ell_A(\sigma a_i)\},
$$
and the {\it A-descent number} of $\sigma\in A_{n+1}$
is defined by
$$
\des_A(\sigma):=|\D_A(\sigma)|.
$$
(note that the strict relation $>$, in the definition of an $S$-descent
in Section~\ref{class},
is replaced in the $A$-analogue by $\ge$).

\item
Define the {\it alternating reverse major index} of $\sigma\in A_{n+1}$ as
$$
\rmaj_{A_{n+1}}(\sigma):=\sum\limits_{i\in \D_A(\sigma)} (n-i) .
$$
\end{enumerate}
\end{df}

\subsubsection{The Delent Number}

New statistics, for the alternating 
group, as well as for the symmetric group, are introduced.

\begin{df} (See Definition~\ref{dfDEL})
\begin{enumerate}
\item
Let $w\in S_n$. The $S$-delent number of $w$ is
the number of times that $s_1=(1,2)$ occurs in the 
$S$ canonical presentation of $w$, and is denoted by $del_S(w)$. 
\item
Let $v\in A_{n+1}$. The $A$-delent number of $v$ is
the number of times that $a^{\pm 1}_1$ 
occur in the $A$ canonical presentation of $v$, and is denoted by $del_A(v)$. 
\end{enumerate}
\end{df}

A combinatorial interpretation of the delent numbers,
$del_S$  and $del_A$, is given in
Section 7. 
Let $w\in S_n$, 
then $j$ is a {\it l.t.r.min}
(left-to-right minimum) of $w$ if $w(i)>w(j)$ for all $1\le i<j$. 

\begin{pro}\label{pro13-i} (see Proposition~\ref{pro8})
For every permutation $w\in S_n$  denote
$$
Del_S(w)= \{1<i\le n|\ i\ \rm{is}\ \rm{a}\ l.t.r.min\}, 
$$ 
then
$$
del_S(w)= | Del_S(w) |.
$$ 
\end{pro}
%
%
Similar to l.t.r.min, we define an
{\it almost left to right minimum} (a.l.t.r.min) of $w\in A_{n+1}$ as follows:

$j$ is an {\it a.l.t.r.min} of $w$ if $w(i)<w(j)$ for 
at most one $j$ less than $i$. Define $Del_A(w)$ as the set
of the {\it almost left-to-right minima} of $w$.
Then $del_A(v)=\mid Del_A(w)\mid$, i.e.~is the number of {\it a.l.t.r.min}
of $w$, see Proposition  7.7.

\medskip


\medskip

%
We also have
\begin{pro}\label{pro12-i}
(See Proposition~\ref{pro12})
Let $w\in A_{n+1}$, then  
$$del_S(w)=\ell_S(w)-\ell_A(w).$$ 
\end{pro}
%

%
%
%


\medskip

\subsubsection{Equi-distribution Identities}

 The covering map $f:A_{n+1} \to S_n$, presented in Definition~\ref{df2},
 allows us to translate $S_n$-identities, which involve the delent statistic, 
 into corresponding $A_{n+1}$-identities.
This strategy is used in the proofs of part (2) of the following theorems.

\medskip

Part (1) of the following theorem is a new generalization 
of MacMahon's classical identity, and part (2) is its $A$-analogue.

\begin{thm} (see Theorem~\ref{thm101})
$$
\sum_{\sg\in S_n} q^{\ell_S(\sg)}t^{del_S(\sg)}=
\sum_{\sg\in S_n} q^{\rmaj_{S_n}(\sg)}t^{del_S(\sg)}= \leqno(1)
$$
$$
=(1+qt)(1+q+q^2t)\cdots (1+q+\ldots  +q^{n-1}t); 
$$
and
$$
\sum_{w\in A_{n+1}}q^{\ell_A(w)}t^{del_A(w )}=
\sum_{w\in A_{n+1}} q^{\rmaj_{A_{n+1}}(w)}t^{del_A(w)}
\leqno(2)
$$
$$
=(1+2qt)(1+q+2q^2t)\cdots (1+q+\ldots + q^{n-2}+2q^{n-1}t ). 
$$
\end{thm}

\bigskip
Recall the standard notation $[m]=\{1,\ldots ,m\}$.
The main theorem in this paper strengthens Theorem~\ref{class1}, 
and also gives its $A$-analogue. This is

\begin{thm}\label{main1}
 (See Theorem 9.1)
For every subsets $D_1\subseteq [n-1]$ and $D_2\subseteq [n]$
$$
\sum\limits_{\{\pi\in S_n|\ \D_S(\pi^{-1})\subseteq D_1,\ 
\Del_S(\pi^{-1})\subseteq D_2\}} q^{\rmaj_{S_n}(\pi)}
=\leqno(1)
$$
$$
\sum\limits_{\{\pi\in S_n|\ \D_S(\pi^{-1})\subseteq D_1,\
\Del_S(\pi^{-1})\subseteq D_2\}} q^{\ell_S(\pi)},
$$
and
$$
\sum\limits_{\{\sigma\in A_{n+1}|\  \D_A(\sigma^{-1})\subseteq D_1,\
\Del_A(\sigma^{-1})\subseteq D_2\}} q^{\rmaj_{A_{n+1}}(\sigma)}
=\leqno(2)
$$
$$
\sum\limits_{\{\sigma\in A_{n+1}|\  \D_A(\sigma^{-1})\subseteq D_1,\
\Del_A(\sigma^{-1})\subseteq D_2\}} q^{\ell_A(\sigma)}.
$$
\end{thm}

\vskip 0.1 truecm

This shows that the delent set and the descent set play a similar role
in these identities.

\bigskip

The  $A$-analogue of Corollary~\ref{class2} follows.
It is obtained as a special case of
Corollary \ref{shuff.2}(2) (by substituting $q_3=1$). 

\begin{cor} (See Corollary~\ref{shuff.2}) 
$$
\sum\limits_{\sigma\in A_{n+1}} q_1^{\rmaj_{A_{n+1}}(\sigma)} q_2^{\des_A(\sigma^{-1})}=
\sum\limits_{\sigma\in A_{n+1}} q_1^{\ell_A(\sigma)}
q_2^{\des_A(\sigma^{-1})}.
$$
\end{cor}

Note that, while the $S$-identity holds for $\maj_{S_n}$ 
as well as for $\rmaj_{S_n}$,
it is not possible to replace $\rmaj_{A_{n+1}}$ by $\maj_{A_{n+1}}$ 
in the $A$-analogue.

\section{Preliminaries}\label{prel}

\subsection{Notation}\label{prel.n}

For an integer $a$ we  let $[a] :=\{ 1,2, \ldots , a \} $ (where $[0]\ :=
\emptyset $).
Let $n_1,\dots,n_{r}$ be non-negative integers such that 
$\sum_{i=1}^{r} n_i=n$.
Recall that the {\it $q$-multinomial coefficient} 
${n \brack n_1,\ldots,n_{r}}_q$ 
is defined by:
\[ 
[0]!_q := 1,
\]
\[
[n]!_q := [n-1]!_q \cdot (1 + q + \ldots + q^{n-1}) \qquad (n \ge 1),
\]
\[
{n \brack n_1 \ldots n_{r}}_q := \frac{[n]!_q}{[n_1]!_q \cdots [n_{r}]!_q}.
\]

Represent $\sg\in S_n$ by `its second row' 
$\sg =[\sg (1),\ldots , \sg (n)]$. 
We also use the cycle-notation; in particular, we denote
$s_i:=(i,i+1)$, the transposition of $i$ and $i+1$. Thus
\begin{eqnarray}\label{eqn3}
[\ldots , \sg (i), \sg (i+1) ,\ldots ]s_i=
[\ldots , \sg (i+1), \sg (i) ,\ldots ]
\end{eqnarray}
(i.e. only $\sg (i), \sg (i+1)$  switch places).

\subsection{The Coxeter System of the Symmetric Group}\label{prel.s}

The symmetric group on $n$ letters, denoted by $S_n$, is generated
by the set of adjacent transpositions $S:=\{(i,i+1)|\ 1\le i<n\}$.

The defining relations of $S$ are the Moore-Coxeter relations :
$$
(s_is_{i+1})^3=1 \qquad (1\le i <n),
$$
$$
(s_is_j)^2=1 \qquad (|i-j|>1)
$$
$$
s_i^2=1 \qquad (\forall i).
$$
This set of generators is called the {\it Coxeter system} of $S_n$.

For $\pi\in S_n$
let $\ell _S(\pi)$ be the standard length of $\pi$ with respect to $S$
(i.e. the length of the canonical presentation of $\pi$,
see Section~\ref{capre} ).
Let $w$ be a word on the letters $S$. 
A {\it commuting move} on $w$ switches the positions of consequent
letters $s_is_j$ where $|i-j|>1$. A {\it braid move} replaces
$s_is_{i+1}s_i$ by $s_{i+1}s_is_{i+1}$ or vice versa.
The following is a well known fact, but we shall not use it in this paper.

\begin{fac}\label{p0}
All irreducible expressions of $\pi\in S_n$ are of length $\ell _S(\pi)$.
For every pair of irreducible words of $\pi\in S_n$,
it is possible to move from one to another along commuting and braid moves.
\end{fac}

\subsection{Permutation Statistics}\label{prel.p}

There are various statistics on the symmetric groups $S_n$, like 
the {\it descent} number and the {\it major} index. We introduce and study 
analogous statistics on the alternating groups $A_n$. As was mentioned, to 
distinguish we add `sub S' and `sub A' accordingly.\\
\newline
Given a permutation $\pi =[\pi(1), \ldots ,\pi(n)]$ in the symmetric group 
$S_n$, we say that a pair $(i,j)$, $1\le i< j\le n$ is an
{\em inversion} of $\pi$ if $\pi(i)>\pi(j)$.
The set of inversions of $\pi$ is denoted by
$Inv_S(\pi)$
and its cardinality is denoted by $\inv_S(\pi)$. Also
$1\le i <n$ is a {\em descent} of $\pi $ if
$\pi(i)>\pi(i+1)$. For the definitions of the descent set $\D_S(\pi)$,
the descent number $\des_S(\pi)$, the major index $\maj_S(\pi)$
and the reverse major index $\rmaj_{S_n}(\pi)$, see Subsection 1.2.
%

\medskip

Note that $i$ is a descent of $\pi$ if and only if 
$\ell _S(\pi s_i)<\ell _S(\pi)$.
Thus (as already mentioned in Subsection 1.2), 
the descent set, and consequently the other related statistics,
have an algebraic interpretation in terms of
the Coxeter system. Also, for every $\pi\in S_n$
\begin{eqnarray}\label{eqn3i}
\inv _S(\pi)=\ell _S(\pi).
\end{eqnarray}

\bigskip

The following well known identity is due to MacMahon
\cite{MM}. See, e.g. \cite{F} and \cite[Corollaries 1.3.10 and 4.5.9]{StI}.

\begin{thm}\label{p1}
$$
\sum\limits_{\pi\in S_n}q^{\inv _S(\pi)}=
\sum\limits_{\pi\in S_n}q^{\maj _S(\pi)}=
$$
$$
=[n]!_q=(1+q)(1+q+q^2)\cdots (1+q+\dots+q^{n-2}+q^{n-1}).
$$
\end{thm}

\smallskip

The following theorem is a reformulation of \cite[Theorem 1]{FS}.

\begin{thm}\label{p2}
For every $B\subseteq [n-1]$,
$$
\sum\limits_{\{\pi\in S_n|\ \D _S(\pi^{-1})=B\}}q^{\inv _S(\pi)}=
\sum\limits_{\{\pi\in S_n|\ \D _S(\pi^{-1})=B\}}q^{\maj _S(\pi)}.
$$
\end{thm}
\medskip
{\bf Note.} Let $\sg\in S_n$, $\sg =[\sg(1),\ldots ,\sg (n)]$. Then
$\sg =[\ldots,k, \ldots ,\ell,\ldots ]$ (i.e. $k$ is left of $\ell$ in
$\sg$ ) if and only if $\sg^{-1} (k) < \sg^{-1} (\ell) $.
\medskip
\vskip 0.15 truecm
\noindent {\bf Shuffles.} Let $1\le i \le n-1$, then  $w\in S_n$ is an 
$\{i\}$-shuffle
if it shuffles $\{1,\ldots , i\}$  with $\{i+1 ,\ldots , n\}$; in other
words, if $1\le a<b \le i$ then $w^{-1}(a) <w^{-1}(b)$, and similarly, if
$i+1\le k<\ell\le n$, then $w^{-1}(k) <w^{-1}(\ell)$.

\noindent{\bf Example.} Let $n=4$ and $B=\{2\}$,
then $\{1,2\}$ and $\{3,4\}$ are being shuffled, hence
$$
[1,2,3,4],\;\; [1,3,2,4],\;\; [1,3,4,2],\;\; [3,1,2,4],\;\;
 [3,1,4,2],\;\; [3,4,1,2]
$$
are all the $\{2\}$-shuffles. 

\vskip 0.15 truecm
\noindent More generally, let $B=\{i_1,\dots,i_k\}\subseteq [n-1]$,
where $i_1<\ldots <i_k$.
Denote $i_0:=0$ and $i_{k+1}:=n$.
A {\it $B$-shuffle} is a permutation which shuffles 
$\{1,\ldots ,i_1\},\{i_1+1,\ldots ,i_2\},\ldots $ 
Thus $\pi\in S_n$ is a $B$-shuffle if it satisfies: if
$i_j\le a<b\le i_{j+1}$ for some $0\le j\le k$, then 
$\pi =[\ldots ,a,\ldots ,b,\ldots]$ (i.e. $a$ is left of $b$ in $\pi$).
Notice that in particular there can be no descent for $\pi^{-1}$ 
on any $a$, $i_{j}<a<i_{j+1}$, hence $Des _S(\pi^{-1})\subseteq B$. The opposite 
is also clear, hence

\smallskip

\begin{fac}\label{p3.1}
For every $B\subseteq [n-1]$
$$
\{\pi\in S_n\ |\ \D _S(\pi^{-1})\subseteq B\}=
\{\pi\in S_n\ |\ \pi\ \hbox{is a $B$-shuffle}\}.
$$
\end{fac}

\medskip

\noindent
For a permutation $\pi\in S_n$ let
$$
\supp(\pi):=\{1\le i\le n\ |\ \pi(i)\not= i \}
$$
be the {\it support} of $\pi$.

\smallskip

Let $k\in [n-1]$, and let
$\pi_1,\pi_2$ be permutations in $S_n$, such that 
$\supp(\pi_1)\subseteq [k]$
and $\supp(\pi_2)\subseteq [k+1,n]$.
A permutation $\sigma\in S_n$ is called a
shuffle of $\pi_1$ and $\pi_2$ 
if $\sigma=\pi_1\pi_2 r$ for some $\{k\}$-shuffle $r$.
Equivalently, $\sigma$ is a shuffle of $\pi_1$ and $\pi_2$ if and only if
the letters of $[k]$ appear in $\sigma$ in the same order as they appear 
in $\pi_1$ 
and the letters of $[k+1,n]$ appear in $\sigma$ in the same order as 
they appear in $\pi_2$.
The following is a special case of \cite[Prop. 1.3.17]{StI}.

\medskip

\begin{fac}\label{p3i}
Let $k\in [n]$, and let
$\pi_1,\pi_2$ be permutations in $S_n$, such that 
$\supp(\pi_1)\subseteq [k]$
and $\supp(\pi_2)\subseteq [k+1,n]$.
Then
$$
\sum\limits_{\D(r^{-1})\subseteq \{k\}} 
q^{\inv _S(\pi_1\pi_2 r)-\inv _S(\pi_1)-\inv _S(\pi_2)}=
{n\brack k}_q. 
$$
\end{fac}

The following analogue is a special case of a well known theorem
of Garsia and Gessel. It should be noted that, while Garsia-Gessel's Theorem
is stated in terms of sequences, our reformulation is in terms of permutations.

\begin{thm}\label{p3}\cite[Theorem 3.1]{GG}
Let $k\in [n-1]$, and let
$\pi_1,\pi_2$ be permutations in $S_n$, such that 
$\supp(\pi_1)\subseteq [k]$
and $\supp(\pi_2)\subseteq [k+1,n]$.
Let 
$\nu_k:= (1,k+1)(2,k+2)\cdots (n-k,n)
\in S_n$.
Then
$$
\sum\limits_{\D _S(r^{-1})\subseteq \{k\}}
q^{\maj _S(\pi_1\pi_2 r)-\maj _S(\pi_1)-\maj _S(\nu_k^{-1}\pi_2\nu_k)}=
{n\brack k}_q. 
$$
\end{thm}

In order to translate Theorem \ref{p3} into Garsia-Gessel's terminology,
note that $\pi_1\pi_2r$ are shuffles of $\pi_1$ and $\pi_2$ (as mentioned 
above); thus the sum runs over all shuffles of $\pi_1$ and $\pi_2$. Also, 
$\maj_S(\nu_k^{-1}\pi_2\nu_k)$ is the major index of $\pi_2$, when it is 
considered as a sequence on the letters $[k+1,n]$.

\medskip

\begin{rem}\label{p4}
In general, it is possible to replace a statement involving $\maj$
by a corresponding statement involving $r\maj$, using
the following automorphism $\sg\to \hat\sg $ :\\
Let $\rho_n\in S_n$ denote the involution
$$
\rho_n:= (1,n)(2,n-1)\cdots  
([n/2 ],[(n+3)/2 ] ),
$$
where for a real number $\alpha$, $[\alpha]$ is the `integer part' of $\alpha$.
Define
$$
\hat\sg := \rho_n\sg\rho_n .
$$
Then $\sg\to \hat\sg$ 
is an automorphism of $S_n$ with the following properties :
\begin{enumerate}
\item
Let $i<j$ and $\sg(i)>\sg(j)$, then $n+1-j<n+1-i$ and 
$\hat\sg(n+1-j)>\hat\sg(n+1-i)$. In particular, $i\in Des _S(\sg)$
if and only if $n-i\in Des _S(\hat\sg)$, hence 
\begin{eqnarray}\label{eq0}
r {maj} _{S_n}(\sg)=maj _{S}(\hat\sg).
\end{eqnarray}

\item
There is a bijection between $Inv _S(\sg)$ and $Inv _S(\hat\sg)$ given by

$(i,j)\leftrightarrow (n+1-j,n+1-i)$, hence 
\begin{eqnarray}\label{eq01}
inv _S(\sg)=inv _S(\hat\sg).
\end{eqnarray}

\item
Part 1 implies that $\sg$ is an $\{i\}$-shuffle if and only if $\hat\sg$ is 
an $\{n-i\}$-shuffle, i.e. 
\begin{eqnarray}\label{eq011}
\D_S(\sg^{-1})\subseteq \{i\}\Longleftrightarrow 
\D_S(\hat\sg^{-1})\subseteq \{n-i\}.
\end{eqnarray}
This easily generalizes to $B$-shuffles.

\end{enumerate}
\end{rem}

\noindent
Note that by (\ref{eq0}) and \cite[Claim 0.4]{R}, for every $\pi\in S_n$
$$
\rmaj_{S_n}(\pi)= charge (\pi^{-1}),
$$
where the charge is defined as in \cite[p. 242]{Md}.

%
%

\section{The $S$ and $A$ Canonical Presentations}\label{capre}

In this section we consider canonical presentations of elements in $S_n$
and in $A_n$ by the corresponding Coxeter generators.
This presentation for $S_n$ is well known,
 see for example~\cite[pp. 61-62]{Go}. The analogous presentation for $A_n$
follows from the properties of the Mitsuhashi's Coxeter generators.

\subsection{The $S_n$ Case}

The $S_n$ canonical presentation is proved below, using the $S$-procedure,
which is also applied later.

\medskip

Recall that
$s_i=(i,i+1)$, $1\le i<n$, are the Coxeter generators of
$S_n$. For each $1\le j\le n-1$ define 
\begin{eqnarray}\label{eqn1}
R^S_j=\{ 1,s_j,s_js_{j-1},\ldots , s_js_{j-1}\cdots s_1\}
\end{eqnarray}
and note that $R^S_1,\ldots ,R^S_{n-1} \subseteq S_n$.
.
\begin{thm}\label{thm1} (see~\cite[pp. 61-62]{Go})
Let $w\in S_n$, then there exist unique elements $w_j\in R^S_j$,
 $1\le j\le n-1$, such that $w=w_1\cdots w_{n-1}$. 
Thus, the presentation $w=w_1\cdots w_{n-1}$ is unique. 
\end{thm}
\begin{df}\label{defCP}
Call the above $w=w_1\cdots w_{n-1}$  in Theorem~\ref{thm1}
the $S$ canonical presentation of $w\in S_n$.
\end{df}
\noindent {\bf A proof} of Theorem~\ref{thm1}
follows from the following S-Procedure.\\ 
\newline
{\bf The S-Procedure.} The following is a simple procedure for calculating the 
$S$ canonical presentation of a given $w\in S_n$. 
It can also be used to prove
Theorem~\ref{thm1}, as well as various other facts. 
%
Let $\sg\in S_n$, $\sg (r)=n$,
$\sg =[\ldots,n,\ldots ]$, then apply Equation~(\ref{eqn3}) to 
`pull $n$ to its place on the right':
$\sg s_rs_{r+1}\cdots s_{n-1}=[\ldots \ldots,n]$. This gives
$w_{n-1}= s_{n-1}\cdots s_{r+1}s_r$. Next, in 
\[
\sg w^{-1}_{n-1}=\sg s_rs_{r+1}\cdots s_{n-1}=
[\ldots,n-1, \ldots,n], 
\]
pull $n-1$ to its right place (second from right) by a similar
product $s_ts_{t+1}\cdots s_{n-2}$. This yields 
$w_{n-2}= s_{n-2} \cdots s_t $. Continue! Finally, $\sg =w_1\cdots w_{n-1}$.

For example, let $\sg =[2,5,4,1,3]$, then $w_{n-1}=w_4=s_4s_3s_2$;
$\sg w^{-1}_4=[2,4,1,3,5]$, therefore $w^{-1}_3=s_2s_3$. Check that 
$w_2=1$ and, finally, $w_1=s_1$. Thus 
$\sg =w_1\cdots w_4=(s_1)(1)(s_3s_2)(s_4s_3s_2)$.
\vskip 0.35 truecm
The uniqueness in Theorem~\ref{thm1} follows by cardinality, since the 
number of canonical
words in $S_n$ is at most 
\[
\prod_{j=1}^{n-1} card(R^S_j)=\mid S_n \mid .
\]
This proves Theorem~\ref{thm1}.
\qed

%
\vskip 0.25 truecm
\subsection{A Generating Set for $A_n$}
We turn now to $A_n$.
As was already mentioned in~\ref{stat}, we let 
$$
a_i:=s_1s_{i+1}\qquad (1\le i\le n-1).
$$
The set 
$$
A:=\{a_i\ | \ 1\le i\le n-1\}
$$
generates the alternating group on $n$ letters $A_{n+1}$.
This generating set and its following properties appear in \cite{Mi}.

\begin{pro}\label{m0}\cite[Proposition 2.5]{Mi}
The defining relations of $A$ are
$$
(a_ia_j)^2=1 \qquad (|i-j|>1) ;
$$
$$
(a_ia_{i+1})^3=1 \qquad (1\le i< n-1) ;
$$
$$
a_1^3=1\qquad\mbox{and}\qquad a_i^2=1 \qquad (1<i\le n-1) .
$$
\end{pro}
The general braid-relation $(a_ia_{i+1})^3=1\;$ implies the following 
braid-relations.
\begin{enumerate}
\item
$a_2a_1a_2=a^{-1}_1a_2a^{-1}_1$ and
\item
$a_2a^{-1}_1a_2=a_1a_2a_1$. 
\item
$a_{i+1}a_{i }a_{i+1}=a_{i }a_{i+1}a_{i }$ if $i\ge 2\;\;$ (since $a^{-1}_i=a_i$).
\end{enumerate}
Let 
$$
\overline A:=A\cup\{a_1^{-1}\},
$$
where $A$ is defined as above. Clearly, $\overline A$ is a generating set 
for $A_{n+1}$. 

\subsection{The $A$ Canonical Presentation}

Mitsuhashi's Coxeter generators are now applied to obtain a unique canonical
presentation for elements in the alternating group.

\medskip

For each $1\le j\le n-1$ define 
\begin{eqnarray}\label{eqn2}
R^A_j=\{ 1,\;a_j,\;a_ja_{j-1},\;\ldots ,\; a_j\cdots a_2,\;a_j
\cdots a_2a_1,\;a_j\cdots a_2a^{-1}_1\}
\end{eqnarray}
and note that $R^A_1,\ldots ,R^A_{n-1} \subseteq A_{n+1}$.

\begin{thm}\label{thm2}
Let $v\in A_{n+1}$, then there exist unique elements $v_j\in R^A_j$,
 $1\le j\le n-1$, such that $v=v_1\cdots v_{n-1}$, and this presentation 
is unique.
\end{thm}
\begin{df}\label{defACP}
Call the above $v=v_1\cdots v_{n-1}$ in Theorem~\ref{thm2}
the $A$ canonical presentation of $v$.
\end{df}
\noindent {\bf Proof} of Theorem~\ref{thm2}.
Let $v=w_1\cdots w_n$, $w_j\in R^S_j$, be the $S$ canonical presentation of $v$. Rewrite that
presentation explicitly as 
\begin{eqnarray}\label{eqn15}
v=(s_{i_1}s_{i_2})\cdots (s_{i_{2r-1}}s_{i_{2r}}).
\end{eqnarray}
Note that $s_is_j=(s_is_1)(s_1 s_j)=a^{-1}_{i-1}a_{j-1}$ (denote $a_0=1$). Thus each $s_i$ in~(\ref{eqn15}) is replaced by
a corresponding $a^{\pm 1}_{i-1}$. It follows that for each $2\le j\le n$,
$w_j$ is replaced by $v_{j-1}\in R^A_{j-1}$ and $v=v_1\cdots v_{n-1}$.
This proves the existence of such a presentation.
\vskip 0.4 truecm
A second proof of the existence follows from the following A-procedure.\\
\newline
{\bf The A-Procedure} is similar to the S-procedure. We describe its
first step, which is also its inductive step. \\
Let $\sg\in A_{n+1}$, $\sg = [\ldots ,n+1,\ldots ]$. As in the S-procedure,
pull $n+1$ to the right:
$\sg s_rs_{r+1}\cdots s_n =[b_1,b_2,\ldots , n+1]$. The (S-) length of
$s_rs_{r+1}\cdots s_n$ is $n-r+1$; if it is odd, use
$\sg s_rs_{r+1}\cdots s_n s_1=[b_2,b_1,\ldots , n+1]$. Thus 
\[
v_{n-1}=\cases
{s_ns_{n-1}\cdots s_r, & \hbox{ if $n-r+1$ is even;}\cr
s_1s_ns_{n-1}\cdots s_r, & \hbox{ if $n-r+1$ is odd.}\cr}
\]
\vskip 0.15 truecm
\noindent {\it The case} $r\ge 2$. Then $s_1s_j=s_js_1$ 
for all $j\ge r+1$, hence
\[
v_{n-1}=\cases
{(s_1s_n)(s_1s_{n-1})\cdots (s_1s_r)=a_{n-1}\cdots a_{r-1}, 
& \hbox{ if $n-r+1$ is even;}\cr
(s_1s_1)(s_1s_n)\cdots (s_1s_r)=a_{n-1}\cdots a_{r-1}, 
& \hbox{ if $n-r+1$ is odd.}\cr}
\]
\noindent {\it The case} $r=1$. If $n-r+1=n$ is even, 
\[
v_{n-1}=s_n\cdots s_2s_1=(s_1s_n)\cdots (s_1s_3)(s_2s_1)
=a_{n-1}\cdots a_2 a^{-1}_1 ,
\]
and similarly if $n-r+1$ is odd.\\
This completes the first step. In the next step, pull $n$ to the $n$-th 
position (i.e. second from the right), etc. This proves the existence of
such a presentation $v=v_1\cdots v_{n-1}$.
\vskip 0.05 truecm
{\it Example.} Let $\sg =[3,5,4,2,1]$, so $n+1=5$. Now 
$\sg s_2s_3s_4=[3,4,2,1,5]$ and since $s_2s_3s_4$ is of odd length (=3), 
permute 3 and 4:  $\sg s_2s_3s_4s_1=[4,3,2,1,5]$. Thus 
$v_3=s_1s_4s_3s_2=(s_1s_1) (s_1s_4)(s_1s_3) (s_1s_2)=a_3a_2a_1$.
Similarly, $v_2=a_2a^{-1}_1$ and $v_1=a_1$, hence
$[3,5,4,2,1]=(a_1)(a_2a^{-1}_1)(a_3a_2a_1)$.
\vskip 0.25 truecm
{\bf Uniqueness} follows by cardinality: note that for all $1\le j\le n-1$,
$\mid R^A_j\mid=j+2$,
hence the number of such words $v_1\cdots v_{n-1}$ in  $A_{n+1}$ is at most
\[
\prod_{j=1}^{n-1}(j+2)=\mid A_{n+1}\mid.
\]
Since each element in $A_{n+1}$ does have such a presentation, 
this implies the 
uniqueness - and the proof of Theorem~\ref{thm2} is complete. \qed
\vskip 0.25 truecm
Given $w\in S_n$, we say that $s_i$ occurs $\ell$ times in $w$ if it
occurs $\ell$ times in the canonical presentation of $w$. Similarly for the number 
of occurrences of $a_i$, or of $a^{-1}_1$, in $v\in A_{n+1}$. 
The number of occurrences of
$s_1$, as well as those of $a^{\pm 1}_1$, are of particular importance 
in this paper. 
\begin{lem}\label{lem9}
Let $w\in S_n$, then the number of occurrences of $s_i$ in $w$ equals 
the number of occurrences of $s_i$ in $w^{-1}$. Similarly for $A_{n+1}$
and $a_1^{\pm 1}$.
\end{lem}
This is an obvious corollary of
\begin{lem}\label{lem2}
Let $w=s_{i_1}\cdots s_{i_p}$ be the  canonical presentation of
$w\in S_n$. Then the canonical presentation of $w^{-1}$ is obtained
from the presentation $w^{-1}=s_{i_p}\cdots s_{i_1}$ by commuting
moves only - without any braid moves.\\
Similarly for $v,\,v^{-1}\in A_{n+1}$.
\end{lem}
{\it Proof.} We prove for $S_n$. The proof is 
by induction on $n$.
Write $w=w_1\cdots w_{n-1}\;$,
$w_j\in R^S_j$. If $w_{n-1}=1$ then $w\in S_{n-1}$ and the proof follows
by induction. 

Let  $w_{n-1}=s_{n-1}s_{n-2}\cdots s_k$ where $1\leq k\leq n-1$. 
Now either $w_{n-2}=1$ or $w_{n-2}=s_{n-2}s_{n-3}\cdots s_{\ell}$ for some 
$1\leq \ell\leq n-2$, and similarly for $w_{n-3}$, $ w_{n-4} $ etc.
The case $w_{n-2}=1$ is similar to the case $w_{n-2}\ne 1$
and is left to the reader, so let $w_{n-2}\ne 1$ and
\[
w^{-1}=w^{-1}_{n-1}w^{-1}_{n-2}\cdots=
(s_k\cdots s_{n-1})(s_{\ell}\cdots s_{n-2})w^{-1}_{n-3}
w^{-1}_{n-4}\cdots
\]
Notice that $s_{n-1}(s_{\ell}\cdots s_{n-3})=
(s_{\ell}\cdots s_{n-3})s_{n-1}$, hence 
\[
w^{-1}=(s_k\cdots s_{n-2})(s_{\ell}\cdots s_{n-3})(s_{n-1}s_{n-2})
w^{-1}_{n-3}w^{-1}_{n-4}\cdots
\]
Next,  move $s_{n-1}s_{n-2}$ to the right, similarly, by commuting moves. 
Continue
by similarly pulling $s_{n-3}$ - in $w^{-1}_{n-3}$ - to the right,
etc.
It follows that by such commuting moves we obtain 
\[
w^{-1}={\overline w\;}^{-1}(s_{n-1}s_{n-2}\cdots s_d)
\]
for some $d$, where $\overline w =s_{j_r}\cdots s_{j_1}\in S_{n-1}$,
and is in canonical form.
By induction, transform ${\overline w\;}^{-1}$ to its canonical form
by commuting moves - and the proof is complete.
\qed\\
\newline
\section{The Lengths Statistics}
The canonical presentations of the previous sections
allow us to introduce the S and the A lengths.
\begin{df}\label{df1}(The length statistics).
\begin{enumerate}
\item
Let $w\in S_n$ with $w=s_{i_1}\cdots s_{i_r}$ its $S$ canonical presentation,
then its S-length is $\ell_S(w)=r$.
\item
Let $v\in A_{n+1}$ with $v=a^{\ep _1}_{i_1}\cdots a^{\ep _r} _{i_r}$
($\ep_i=\pm 1$) its $A$ canonical presentation,
then its A-length is $\ell_A(v)=r$.
\end{enumerate}
\end{df}
For example, $\ell_A(a_1)=1$ and  $\ell_S(a_1)=\ell_S(s_1s_2)=2$.

\begin{rem}\label{r0}
An analogue of Fact \ref{p0} holds :
All irreducible expressions of $v\in A_{n-1}$ are of length $\ell _A(v)$.
This fact will not be used in the paper.
\end{rem}
%
\begin{df}\label{dfDEL}
\begin{enumerate}
\item
Let $w\in S_n$. The number of times that $s_1$ occurs in the $S$ canonical presentation
of $w$ is denoted by $del_S(w)$. 
\item
Let $v\in A_{n-1}$. The number of times that $a^{\pm 1}_1$ 
occurs in the $A$ canonical presentation of $v$ is denoted by $del_A(v)$. 
\end{enumerate}
\end{df}
For example, $del_S(s_1s_2s_1s_3)=2$ and 
$del_A(a^{-1}_1a_2a_1a_3a_2a^{-1}_1))=3$.\\
\newline
A combinatorial characterization of $del_S$
($del_A$) is given in section~\ref{del}.

\smallskip

\noindent
Relations between $del_S$ and the S and the A lengths of $v\in A_{n+1}$
are given by the following proposition.
\begin{pro}\label{pro12}
Let $w\in A_{n+1}$, then  
$$\ell_A(w)=\ell_S(w)-del_S(w).$$ 
Moreover, let 
\begin{eqnarray}\label{eqn12}
w=s_{i_1}\cdots s_{i_{2r}}=w_1\cdots w_n,\qquad w_i\in R^S_i,
\end{eqnarray} be its $S$ canonical presentation and 
\begin{eqnarray}\label{eqn13}
w=a^{\ep _1}_{j_1}\cdots a^{\ep _t}_{j_t} =  v_1\cdots v_{n-1},
\qquad v_i\in R^A_i,
\end{eqnarray}
its $A$ canonical presentation. Then
\begin{eqnarray}\label{eqn14}
\ell _A(v_i)=\cases
{\ell _S(w_{i+1})  & \hbox{ if $s_1$ does not occur in $w_{i+1}$;}\cr
\ell _S(w_{i+1})-1  & \hbox{ if $s_1$  occurs in $w_{i+1}$.}\cr}
\end{eqnarray}
\end{pro}
%
%
%
{\it Proof.} As in the proof of Theorem~\ref{thm2}, 
the proof easily follows from~(\ref{eqn15}) by replacing
$s_is_j$ by $(s_is_1)(s_1s_j)$.

\qed
%


%
%
The S-lengths $\ell _S(w_{i+1})$ and the A--lengths $\ell _A(v_i)$
in~(\ref{eqn14}) can be calculated directly from 
$w=[b_1,\ldots ,b_{n+1}]$ as follows.
\begin{pro}\label{pro13}
Let $w\in S_{n+1}$ as above.
For each $2\le j\le n$ let $T_j(w)$ denote the set of indices $i$
such that  $i<j$ and
$w=[\ldots,j,\ldots ,i,\ldots ]$ (i.e. $w^{-1}(i)>w^{-1}(j)$); denote
$t_j(w)=\mid T_j(w) \mid$. Keeping the notations of Proposition~\ref{pro12}
we have:
\begin{enumerate}
\item
$\ell _S(w_{j})=t_{j+1}(w)$. Moreover,
$T_{ j+1}(w)$ is the full set $\{1,\ldots ,j\}$ (i.e.
$t_{j+1}(w)=j$) if and only if $s_1$ occurs in $w_{j}$. 
\item
$\ell _A (v_k)$ equals  $\mid T_k(w) \mid$,
provided that $ T_k(w) $ is not the full set $\{1,\ldots ,k-1\}$,
and it equals $\mid T_k(w) \mid -1$ otherwise.
\end{enumerate}
\end{pro}
{\it Proof.} By an easy induction on $n$, prove that 
\[
(\ell _S(w_1),\ldots ,\ell _S(w_n))=(t_2(w),\ldots , t_{n+1}(w)).
\]
This follows since \\
$[b_1,\ldots , b_n,n+1]s_ns_{n-1}\cdots s_r=
[b_1,\ldots , b_{r-1},n+1,b_r,\ldots , b_n]$.\\
Here are the details: 
Write $w=w_1\cdots w_n$, let
$\sg = w_1\cdots w_{n-1}$, so
$\sg = [d_1,\ldots ,d_n,n+1]$. If $w_n=1$, the claim follows by induction.
Let $w_n=s_ns_{n-1}\cdots s_r$ for some $r\ge 1$. Then 
$w=\sg w_n=[d_1,\ldots ,d_{r-1},n+1,d_r,\ldots d_n]$. Thus
$t_{n+1}(w)=n-r+1=\ell _S(w_n)$. Also, for $2\le j \le n$, 
$t_j(w)=t_j(\sg)$, and the proof of part 1 follows by induction.
Part 2 now follows from~(\ref{eqn14}).
\qed 
\section{$f$-Pairs of Statistics}

\subsection{The Covering Map}

Theorems~\ref{thm1} and~\ref{thm2} allow us to introduce
the following definition.
\begin{df}\label{df2} 
Define  $f:A_{n+1}\to S_n$ as follows. 
\[
f(a_1)=f(a^{-1}_1)=s_1\qquad\mbox{and}\qquad f(a_i)=s_i,
\qquad 2\le i\le n-1.
\]
Now extend $f:R^A_j\to R^S_j$ via
\[
f(a_ja_{j-1}\cdots a_{\ell})=s_js_{j-1}\cdots s_{\ell},\qquad
f(a_j\cdots a_1)=f(a_j\cdots a^{-1}_1)=s_j\cdots s_1.
\]
Finally, let $v\in A_{n+1},\quad v=v_1\cdots v_{n-1}$ its unique
$A$ canonical presentation, then 
\[
f(v)=f(v_1)\cdots f(v_{n-1})\,
\]
which is clearly the $S$ canonical presentation of $f(v)$.
\end{df}
Notice that for $v\in A_{n+1}$, $\ell_A(v)=\ell_S(f(v))$. We therefore
say that the pair of 
the length statistics $(\ell_S,\ell_A)$ is an $f$-pair. More generally, we have
\begin{df}\label{df3}
Let $m_S$ be a statistic on the symmetric groups and 
$m_A$ a statistic on the alternating groups. We
say that $(m_S,m_A)$ is an
$f$-pair (of  statistics) if for any $n$ and 
$v\in A_{n+1}$, $m_A(v)=m_S(f(v))$.
\end{df}
Examples of $f$-pairs are given in Proposition~\ref{pair2} below.

\begin{pro}\label{pair1}
Recall Definition~\ref{df8}.1.
For every $\pi\in A_{n+1}$ 
$$
\D_A(\pi)=Des_S(f(\pi)).
$$
\end{pro}
\noindent{\bf Proof -} is left to the reader.
\qed

\medskip

It follows that the descent statistics are $f$-pairs.
By Definition \ref{dfDEL}, $(del_S,del_A)$ is an $f$-pair. We summarize:

\begin{pro}\label{pair2}
The following pairs
\[ (\ell_S,\ell_A), \]
\[ (\des_S,\des_A), \]
\[ (\maj_S,\maj_A), \]
\[ (\rmaj_{S_n},\rmaj_{A_{n+1}}) \]
and
\[ (del_S,del_A) \]
are $f$-pairs.
\end{pro}

\subsection{The `del' Statistics}

The following basic properties of $del_S$
play an 
important role in this paper.

\begin{pro}\label{exa1} 
\begin{itemize}
\item[1.]
For each $w\in S_n\;$, $\mid f^{-1}(w)\mid = 2^{del_S(w)}$. 
\item[2.] For each $w\in S_n$ and $v\in A_{n+1}$
\begin{eqnarray}\label{eqn9}
del_S(w)=del_S(w^{-1})\qquad\mbox {and}\qquad del_A(v)=del_A(v^{-1}).
\end{eqnarray}
\end{itemize}
\end{pro}
\noindent{\bf Proof.} Part 1 follows
since each occurrence of $s_1$ can be replaced by an occurrence of either 
$a_1$ or  $a^{-1}_1$.
Part 2 follows from Lemma~\ref{lem9}.

\qed 

We have the following general proposition.
\begin{pro}\label{pro1}
Let $(m_S,m_A)$ be an $f$-pair of statistics, then for all $n$
\[
\sum_{v\in A_{n+1}}q^{m_A(v)}t^{del_A(v)}=
\sum_{w\in S_{n}}q^{m_S(w)}(2t)^{del_S(w)}.
\]
\end{pro}
{\it Proof.} 
Since $A_{n+1}=\cup_{w\in S_n} f^{-1}(w)$, a disjoint union, we have:
\[
\sum_{v\in A_{n+1}}q^{m_A(v)}t^{del_A(v)}=
\sum_{w\in S_{n}}\sum_{v\in f^{-1}(w)}q^{m_A(v)}t^{del_A(v)}=  
\]
\[
\sum_{w\in S_{n}}\sum_{v\in f^{-1}(w)}q^{m_S(f(v))}t^{del_S(f(v))}=
\sum_{w\in S_{n}}\sum_{v\in f^{-1}(w)}q^{m_S(w)}t^{del_S(w))}=
\]
\[
\sum_{w\in S_{n}}2^{del_S(w)}q^{m_S(w)}t^{del_S(w))}.
\]
\qed\\
A refinement of Proposition~\ref{pro1} is given in Proposition~\ref{pro5}

\begin{pro}\label{pro2}
With the above notations we have:
\begin{enumerate}
\item
$$
\sum_{\sg\in S_n }q^{\ell_S(\sg )}t^{del_S(\sg )}
=(1+qt)(1+q+q^2t)\cdots (1+q+\ldots  +q^{n-1}t).
$$
\item
$$
\sum_{w\in A_{n+1}}q^{\ell_A(w)}t^{del_A(w )}
=(1+2qt)(1+q+2q^2t)\cdots (1+q+\ldots + q^{n-2}+2q^{n-1}t ).
$$
\end{enumerate}
\end{pro}
{\it Proof.}
\begin{enumerate}
\item
The proof of part 1 is similar to the proof of Corollary 1.3.10 
in~\cite{StI}.
Let $w_j\in R^S_j$, then $del_S(w_j)=1$ if $w_j=s_j\ldots s_1$ and $=0$
otherwise. Let $w\in S_n$ and let $w=w_1\cdots w_{n-1}$ be its
$S$ canonical presentation, then
$del_S(w)=del_S(w_1)+\ldots + del_S(w_{n-1})$ and
$\ell_S(w)=\ell_S(w_1)+\cdots + \ell_S(w_{n-1})$. Thus
\[
\sum_{w\in S_n}q^{\ell_S(w)}t^{del_S(w)}=
\prod _{j=1}^{n-1}\left ( \sum _{w_j\in R^S_j}q^{\ell_S(w_j)}t^{del_S(w_j)}\right ).
\]
The proof now follows since 
\[
\sum _{w_j\in R^S_j}q^{\ell_S(w_j)}t^{del_S(w_j)}=1+q+q^2+\ldots +q^{j-2}+q^{j-1}t.
\]
\item
By Proposition~\ref{pro1}, part 2 follows from part 1.
\end{enumerate}\qed
\vskip 0.25 truecm

\subsection{Connection with the Stirling Numbers}

Recall that $c(n,k)$ is the number of permutations in $S_n$ with 
exactly $k$ cycles,
$1\leq k\leq n$: $c(n,k)$ are {\it the sign-less Stirling numbers 
of the first kind}.
Let $w_S(n,\ell)$ denote the number of $S$ canonical words in $S_n$ with 
$\ell$ appearances of $s_1$. 
Similarly, let $w_A(n+1,\ell)$ denote the number of $A$ canonical words 
in $A_{n+1}$ 
with $\ell$ appearances of $a^{\pm 1}_1$. 

\noindent We prove 
\begin{pro}\label {pro3}
Let $0\leq \ell\leq n-1$, then 
\begin{enumerate}
\item
\[
\sum _{\ell\geq 0}w_S(n,\ell)t^\ell=(t+1)(t+2)\cdots (t+n-1),
\]
hence $w_S(n,\ell)=c(n,\ell +1)$.
\item
\[
\sum _{\ell\geq 0}w_A(n,\ell)t^\ell=(2t+1)(2t+2)\cdots (2t+n-1),
\]
hence $w_A(n+1,\ell)=2^\ell \cdot c(n,\ell +1)$.
\end{enumerate}
\end{pro}
{\it Proof.} Substitute $q=1$ in Proposition~\ref{pro2} and,
in part 1, apply 
Proposition 1.3.4 of~\cite{StI}, which states that
\[
\sum_{k=0}^nc(n,k)x^k=x(x+1)(x+2)(\cdots (x+n-1).
\]
\qed 
\vskip 0.25 truecm
further connections with the Stirling numbers are given below
(Propositions~\ref{pro6}, \ref{pro7} and~\ref{pro10}) and in~\cite{RR}.

\subsection{A Multivariate Refinement}

\begin{df}\label{df5}
Let $w\in S_n$, $w=w_1\cdots w_{n-1}$ its $S$ canonical presentation and
let $1\le j \le n-1$. Denote $\ep_{_{S,j}}(w)=1$ if $s_1$ occurs in 
$w_j$, and $\;=0$ otherwise; also denote 
$${\bar \ep}_S(w)=(\ep _{S,1}(w),\ldots ,\ep _{S,n-1}(w))$$
and
$$
t^{{\bar \ep}_S(w)}=t_1^{\ep _{S,1}(w)}\cdots t_{n-1}^{\ep _{S,n-1}(w)}.
$$
Similarly for $v=v_1\cdots v_{n-1}\in A_{n+1}$:
$\ep_{_{A,j}}(v)=1$ if $a^{\pm 1}_1$ occurs in $v_j$, and $\;=0$ otherwise,
and define ${\bar \ep}_A(w)$ similarly. 
%
%
Clearly, $del_S(w)=\sum_j \ep_{_{S,j}}(w)$ and 
$del_A(v)=\sum_j \ep_{_{A,j}}(v)$. 
\end{df}
Proposition~\ref{pro1} admits the following generalization.
\begin{pro}\label{pro5}
Let $(m_S,m_A)$ be an $f$-pair of statistics, then for all $n$
\[
\sum_{v\in A_{n+1}}q^{m_A(v)}\prod _{j=1}^{n-1}
t_j^{\ep_{_{A,j}}(v)}=
\sum_{w\in S_{n}}q^{m_S(w)}\prod _{j=1}^{n-1} 
(2t_j)^{\ep_{_{S,j}}(w)}.
\]
\end{pro}
{\bf The proof} is a slight generalization of the proof of 
Proposition~\ref{pro1} - and is left to the reader.
\vskip 0.25 truecm
We end this section with another two multivariate generalizations,
which will not be used in the rest of the paper.
Proposition~\ref{pro2} generalizes as follows.
\begin{pro}\label{pro6}
Let $\ell_S,\; \ell_A$ be the length statistics, then 
\begin{enumerate}
\item
\[
\sum_{w\in S_{n}}q^{\ell_S(w)}\prod _{j=1}^{n-1}
(t_j)^{\ep_{_{S,j}}(w)}=
(1+qt_1)(1+q+q^2t_2)\cdots (1+q+\ldots  +q^{n-1}t_{n-1}).
\]
\item
\[
\sum_{v\in A_{n+1}}q^{\ell_A(v)}\prod _{j=1}^{n-1}
t_j^{\ep_{_{A,j}}(v)}=
(1+2qt_1)
\cdots (1+q+\ldots + q^{n-2}+2q^{n-1}t_{n-1} ).
\]
\end{enumerate}
\end{pro}
One can generalize Proposition~\ref{pro3}
as follows. Let $w=w_1\cdots w_{n-1}\in S_n$,
a canonical presentation, with $\ep_{_{S,j}}(w)$ and
$\;\overline \ep_{_S}(w)$ as in Definition~\ref{df5}.
Given 
$\overline \ep =(\ep_1,\cdots , \ep_{n-1})$ with all $\ep_i\in\{0,1\}$, 
denote $w_S(n,\overline \ep )=
card \{w\in S_n\mid \overline \ep_S(w)=\overline \ep\}$. Also denote 
$\mid \overline \ep\mid =\sum_j \ep _j$ and $t^{\overline \ep}=
\prod_jt_j^{\ep _j}$. Note that 
\[
\sum_{\mid \overline \ep\mid =\ell} w_S(n,\overline \ep)=
w_\ell(n,\ell )=c(n,\ell +1).
\]
Similarly, introduce the analogous notations for $A_{n+1}$.\\
Proposition~\ref{pro3} now generalizes as follows.
\begin{pro}\label{pro7} With the above notations
\begin{enumerate}
\item
\[
\sum_{\overline \ep}w_S(n,\overline \ep)\,t^{\overline \ep}=
(t_1+1)\cdots (t_{n-1}+n-1).
\]
\item
\[
\sum_{\overline \ep}w_A(n,\overline \ep)\,t^{\overline \ep}=
(2t_1+1)\cdots (2t_{n-1}+n-1).
\]
\end{enumerate}
\end{pro}

\section{The Major Index and the Delent Number }
Recall the definitions of $\rmaj_{S_n}$ and $\rmaj_{A_{n+1}}$
from Subsections 1.2 and 1.3.
In this section we prove
%

\begin{thm}\label{thm101}
$$
\sum_{\sg\in S_n} q^{\ell_S(\sg)}t^{del_S(\sg)}=
\sum_{\sg\in S_n} q^{\rmaj_{S_n}(\sg)}t^{del_S(\sg)}= \leqno(1)
$$
$$
=(1+qt)(1+q+q^2t)\cdots (1+q+\ldots  +q^{n-1}t); 
$$
and
$$
\sum_{w\in A_{n+1}}q^{\ell_A(w)}t^{del_A(w )}=
\sum_{w\in A_{n+1}} q^{\rmaj_{A_{n+1}}(w)}t^{del_A(w)}
\leqno(2)
$$
$$
=(1+2qt)(1+q+2q^2t)\cdots (1+q+\ldots + q^{n-2}+2q^{n-1}t ). 
$$
\end{thm}

Note that Theorem~\ref{thm101} follows from our main theorem~\ref{shuff.1}.
However, the proof of Theorem~\ref{shuff.1} applies the machinery 
required for the proof of Theorem~\ref{thm101} combined with additional, more 
elaborate arguments - therefore we prove it here.

\medskip

Comparing the coefficients of $t^k$ in both parts,
we obtain 
\begin{thm}\label{thm103}
Let $B^S_{n,k} :=\{\sg\in S_n \mid del_S(\sg )=k\}$ and\\
$B^A_{n+1,k} :=\{\sg\in A_{n+1} \mid del_A(\sg )=k\}$.
Then
for each $0\leq k\leq n-1$,
$$
\sum_{\sg\in B^S_{n,k} } q^{\ell_S(\sg)}=
\sum_{\sg\in B^S_{n,k}} q^{\rmaj_{S_n}(\sg)};\leqno(1)
$$
and
$$
\sum_{\sg\in B^A_{n+1,k} } q^{\ell_A(\sg)}=
\sum_{\sg\in B^A_{n,k}} q^{\rmaj_{A_{n+1}}(\sg)}.\leqno(2)
$$
\end{thm}
Note that part 1 is a refinement of MacMahon's equi-distribution theorem.

\medskip


%
%
%
%
The proof of Theorem~\ref{thm101}
follows from the lemmas below. Recall that the descent set - hence 
also the major-indices $maj_S$ and ${\rmaj}_{S_n}$  - are defined for 
any sequence of integers, not necessarily distinct. 
Here $n$ denotes the number of letters in the sequence.

\begin{lem}\label{lem102}
Let $x_1,\ldots , x_n$ and $y$ be integers, not necessarily distinct,
such that $x_i < y$ for $1\leq i \leq n$. Let $u$ be the $n$-tuple 
$u=[x_1,\ldots , x_n]$, and let
\[
v_i=[x_1,\ldots ,x_{i-1},y,x_{i},\ldots , x_n],\qquad 1\leq i \leq n+1
\]
(thus
$v_1=[y,x_1,\ldots , x_n]\quad$ and $\quad v_{n+1}=[x_1,\ldots , x_n,y]$).
Then 
\begin{enumerate}
\item
\begin{eqnarray}\label{eqn4}
\sum_{i=1}^{n+1} q^{maj_S(v_i)}
=q^{maj_S(u)}(1+q+\ldots + q^n)
\end{eqnarray}
\begin{eqnarray}\label{eqn5}
\mbox{and} \qquad\qquad \sum_{i=1}^{n} q^{maj_S(v_i)}
=q^{maj_S(u)}(q+q^2+\ldots + q^n).
\end{eqnarray}
\item
\begin{eqnarray}\label{eqn6}
\sum_{i=1}^{n+1} q^{{\rmaj_{S_{n+1}}}(v_i)}
=q^{{\rmaj_{S_{n}}}(u)}(1+q+\ldots + q^n)
\end{eqnarray}
\begin{eqnarray}\label{eqn7}
\mbox{and}\qquad\qquad  \sum_{i=2}^{n+1} q^{{\rmaj_{S_{n+1}}}(v_i)}
=q^{{\rmaj_{S_{n}}}(u)}(1+q+\ldots + q^{n-1}).
\end{eqnarray}
\end{enumerate}
\end{lem}
Part 1 of Lemma \ref{lem102} is well known. The proof of part 2 is similar.
For the sake of completeness the proof is included.

\medskip

\noindent
{\it Proof.} Denote $u'=[x_1,\ldots , x_{n-1}]\;$ and
$u''=[x_2,\ldots , x_{n}]$. Similarly, denote
$v'_i=[x_1,\ldots ,x_{i-1},y,x_{i},\ldots , 
x_{n-1}],\quad 1\le i\le n,\;$  and\\
$v''_i=[x_2,\ldots ,x_{i-1},y,x_{i},\ldots , x_{n}],\quad 2\le i\le n+1$ \\
(thus $v'_{n}=[x_1,\ldots , x_{n-1},y]$ and 
$v''_{2}=[y,x_2,\ldots , x_{n}]$).
\begin{enumerate}
\item
We argue by induction on $n$, proving~(\ref{eqn4}) and~(\ref{eqn5})  first.
\begin{enumerate}
\item
Assume $x_{n-1}\leq x_n$, then $maj_S(v_{n+1})=maj_S(u)=maj_S(u')=maj_S(v'_{n})$,
$maj_S(v_{n})=maj_S(u)+n$ and $maj_S(v_i)=maj_S(v'_i)$ for $1\leq i\leq n-1$.
It follows that 
\begin{eqnarray*}
\lefteqn
{\sum_{i=1}^{n+1} q^{maj_S(v_i)}=
q^{maj_S(v_{n+1})}+q^{maj_S(v_{n})}+\sum_{i=1}^{n-1} q^{maj_S(v_i)}=} \\
& &  q^{maj_S(u)}+q^nq^{maj_S(u)}+\sum_{i=1}^{n-1} q^{maj_S(v'_i)}=\\
& & q^nq^{maj_S(u)}+\sum_{i=1}^{n} q^{maj_S(v'_i)} =\quad \mbox{(by induction)}\\
& & q^{maj_S(u)}(1+q+\ldots +q^{n-1})+q^nq^{maj_S(u)}.
\end{eqnarray*}
\item
Assume $x_{n-1}> x_n$, then $maj_S(v_{n+1})=maj_S(u)=maj_S(u')+n-1$ and
$maj_S(v_i)=maj_S(v'_i)+n$ for $1\leq i\leq n$.
Thus
\begin{eqnarray*}
\lefteqn{\sum_{i=1}^{n+1} q^{maj_S(v_i)}=q^{maj_S(u)}+
\sum_{i=1}^{n}q^{n+maj_S(v'_i)}=}\\
& &  q^{maj_S(u)}+q^{n}\sum_{i=1}^{n}q^{maj_S(v'_i)}=\quad\mbox{(by induction)}\\
& & q^{maj_S(u)}+q^{n}q^{maj_S(u')}(1+q+\ldots +q^{n-1})=\\
& & q^{maj_S(u)}+q^{maj_S(u)+1}(1+q+\ldots +q^{n-1})=\\
& & q^{maj_S(u)}(1+q+\ldots +q^{n}). 
\end{eqnarray*}
Together, (a) and (b) prove Equation~(\ref{eqn4}). 
Note that Equation~(\ref{eqn5}) follows
from Equation~(\ref{eqn4}) since, in both cases above, 
$maj_S(v_{n+1})=maj_S(u)$.
\end{enumerate}
%
\item
We prove now~(\ref{eqn6}) and~(\ref{eqn7}).
\begin{enumerate}
\item
Assume $x_1\leq x_2$.
First, note that $1\in Des_S(v_1)$ and it contributes
$n+1-1=n$ to $\rmaj_{S_{n+1}}(v_1)$. Let $2\leq k \leq n-1$, then 
$k\in Des_S (u)$ if and only if $k+1\in Des_S (v_1)$; such $k$ contributes 
$n-k=(n+1)-(k+1)$ to both $\rmaj_{S_{n+1}}(v_1)$ and to $\rmaj_{S_{n}}(u)$.
It follows that  $\rmaj_{S_{n+1}}(v_1)=\rmaj_{S_{n}}(u)+n$. 

By similar arguments $\rmaj_{S_{n+1}}(v_i)=\rmaj_{S_{n}}(v''_i)$
for $\;2\leq i \leq n+1$ and also 
$\rmaj_{S_{n}}(u)=\rmaj_{S_{n-1}}(u'')$.
Thus 
\begin{eqnarray*}
\sum_{i=1}^{n+1} q^{\rmaj_{S_{n+1}}}{(v_i)}=q^{\rmaj_{S_{n}}(u)+n}+
\sum_{i=2}^{n+1}q^{\rmaj_{S_{n}}(v''_i)\quad\mbox{(by induction)} }\\
=q^{\rmaj_{S_{n}}(u)+n}+q^{\rmaj_{S_{n-1}}(u'')}
(1+q+\ldots +q^{n-1})
\end{eqnarray*}
and the proof follows.
\item
Assume $x_1 > x_2$. Here (again) $\rmaj_{S_{n+1}}(v_1)=\rmaj_{S_{n}}(u)+n$
while $\rmaj_{S_{n+1}}(v_2)=\rmaj_{S_{n}}(v''_2)=\rmaj_{S_{n}}(u)$.
By similar arguments as above, 
$\rmaj_{S_{n+1}}(v_i)=\rmaj_{S_{n}}(v''_i)+n$
for $3\leq i \leq n+1$.
Also $\rmaj_{S_{n}}(u)=\rmaj_{S_{n-1}}(u'') +n-1$.
Thus 
\begin{eqnarray*}
\lefteqn{\sum_{i=1}^{n+1} q^{\rmaj_{S_{n+1}}(v_i)}=q^{\rmaj_{S_{n}}(u)+n}+
q^{\rmaj_{S_{n}}(u)}+
q^n\sum_{i=3}^{n+1}q^{\rmaj_{S_{n}}(v''_i)}=}\\
& &  
q^{\rmaj_{S_{n}}(u)}+
q^n\sum_{j=2}^{n+1}q^{\rmaj_{S_{n}}(v''_j)}=\mbox{(by induction)}\\
& & 
= q^{\rmaj_{S_{n}}(u)}+
q^{\rmaj_{S_{n-1}}(u'')+n}
(1+q+\ldots +q^{n-1})=\\
& & =q^{\rmaj_{S_{n}}(u)}+q^{\rmaj_{S_{n}}(u)+1}(1+q+\ldots +q^{n-1}),
\end{eqnarray*}
which implies the proof.
Together, (a) and (b) prove Equation~(\ref{eqn6}). 
Now Equation~(\ref{eqn7}) follows
from Equation~(\ref{eqn6}) since, in both cases (a) and (b) above, 
$\rmaj_{S_{n+1}}(v_1)=\rmaj_{S_{n}}(u)+n$.
\end{enumerate}
\end{enumerate}
\qed
%
%
%
\begin{lem}\label{lem103} 
Recall that $R^S_n=\{1,s_n,\ldots ,s_ns_{n-1}\cdots s_1\}\subseteq S_{n+1}$ 
and let $w\in S_n$
(so $w\in S_{n+1}$, where $w(n+1)=n+1$). Then

\[
\sum_{\tau\in R^S_n}q^{maj_S (w\tau)}=
q^{maj _S(w)}(1+q+\ldots + q^n),
\]
and 
\[
\sum_{\tau\in R^S_n}q^{\rmaj_{S_{n+1}} (w\tau)}=
q^{\rmaj_{S_{n}} (w)}(1+q+\ldots + q^n).
\]
\end{lem}
{\it Proof.} 
Write $w\in S_n$ as $w=[w(1),\ldots ,w(n)]\;(=u$ in~\ref{lem102}). 
Similarly write 
$w\in S_n\subseteq S_{n+1}$ as $w=[w(1),\ldots ,w(n),n+1]\;(=v_{n+1}$,
in~\ref{lem102}, where $y=n+1$).  Thus
$$ws_n=[w(1),\ldots ,n+1,w(n)]\quad (=v_n),$$ 
$$ws_ns_{n-1}=[w(1),\ldots ,n+1,w(n-1),w(n)]\quad (=v_{n-1}),$$
etc, and the proof follows by the previous lemma. \qed\\
\newline
\begin{rem}\label{rem2}
Let ${\tilde R}^S_n= R^S_n-\{s_ns_{n-1}\cdots s_1\}\subseteq S_{n+1}$, and let
$\sg\in S_n$. It follows from
Equation~(\ref{eqn7})
that 
\[
\sum_{\tau\in{\tilde R}^S_n}q^{\rmaj_{S_{n+1}} (\sg\tau)}=
q^{\rmaj_{S_{n}} (\sg)}(1+q+\ldots + q^{n-1}).
\]
\end{rem}
\vskip 0.45 truecm
\begin{lem} \label{lem104}
For every $\sg\in S_n$
\[
\sum_{\tau\in R^S_n}q^{\rmaj_{S_{n+1}} (\sg\tau)}t^{del_S(\sg\tau)} =
q^{\rmaj_{S_{n}} (\sg)}t^{del_S(\sg)}(1+q+\ldots + q^{n-1}+ tq^n).
\]
\end{lem}
{\it Proof.} 
By Lemma~\ref{lem103} 
\[
\{\rmaj_{S_{n+1}}(\sg\tau)\mid \tau\in R^S_n\}=
\{\rmaj_{S_{n}}(\sg )+ i \mid 0\leq i \leq n\}.
\]
Let $\eta = s_ns_{n-1}\cdots s_1$ and note that 
$\rmaj_{S_{n}}(\sg )+ n=\rmaj_{S_{n+1}}(\sg\eta )$ (this is the
statement ``$\rmaj_{S_{n+1}}(v_1)=\rmaj_{S_{n}}(u)+n$'' in the proof
of Lemma~\ref{lem102}).

Let $\tau\in R^S_n$. 

If $\tau\ne \eta$ then 
$del_S(\sg\tau)=del_S(\sg)$ since both $ \sg$ and $\sg\tau$ have the same
number of occurrences of $s_1$. By a similar reason
$del_S(\sg\eta)=del_S(\sg)+1$. Thus 
\begin{eqnarray*}
\lefteqn{\{\rmaj_{S_{n+1}}(\sg\tau)del_S(\sg\tau)\mid \tau\in R^S_n\}=}\\
& & \{\rmaj_{S_{n+1}}(\sg\tau)del_S(\sg\tau)\mid \tau\in R^S_n,\;\tau\ne\eta\}
\cup\{\rmaj_{S_{n+1}}(\sg\eta )del_S(\sg\eta )\}=\\
& & \{(\rmaj_{S_{n}}(\sg) +i)del_S(\sg )\mid 0\leq i\leq n-1\}
\cup\{(\rmaj_{S_{n}}(\sg) +n)(del_S(\sg )+1)\}
\end{eqnarray*}
(disjoint unions with no repetitions in the sets)
which translates to 
\[
\sum_{\tau\in R^S_n}q^{\rmaj_{S_{n+1}} (\sg\tau)}t^{del_S(\sg\tau)} =
q^{\rmaj_{S_{n}} (\sg)}t^{del_S(\sg)}(1+q+\ldots + q^{n-1}+ tq^n).
\]   \qed  
\vskip 0.35 truecm
\begin{pro}\label{prop103}
For all $n$
\[
\sum_{\sg\in S_n} q^{\rmaj_{S_n}(\sg)}t^{del_S(\sg)}=
(1+tq)(1+q+tq^2)\cdots (1+q+\ldots + q^{n-2}+tq^{n-1}).
\]
\end{pro}
{\it Proof.} Follows from Lemma~\ref{lem104} by induction on $n$, since
\[
S_{n+1}=\cup_{\tau\in R^S_n}S_n \tau .
\]  \qed\\
\newline
{\bf The proof of Theorem~\ref{thm101}}.\\
Part (1) clearly follows by comparing  
part 1 of Proposition~\ref{pro2} with Proposition~\ref{prop103}.\\
Part (2) follows from part (1) by Proposition~\ref{pro1}. \qed
%
%
\vskip 0.45 truecm

\section{Additional Properties of the Delent Number}\label{del}
We show first that $del_S(w)$ is the number of left-to-right minima of $w$.
\begin{df}\label{df4}
Let $w\in S_n$, 
then $j$ is a l.t.r.min
(left-to-right minimum) of $w$ if $w(i)>w(j)$ for all $1\le i<j$. 
Write $w =[b_1,\ldots ,b_n]$, then $i=1$ is a l.t.r.min
and so is $i$ such that $b_i=1$. We slightly modify the definition, so that
the identity $e=[1,\ldots ,n]$ has no l.t.r.min. This can be done in one of 
the following two ways.\\
Either:\\
1. Do not count $i=1$  as a l.t.r.min (which is Definition~\ref{df4}.1
of l.t.r.min),\\
or:\\
2. Do not count $i$ such that $b_i=1$ as a l.t.r.min
(which is Definition~\ref{df4}.2
of l.t.r.min).\\
3. Define $Del_S(w)$ as the l.t.r.min according to Definition~\ref{df4}.1:
\[
Del_S(w):=\{1< i\le n\ |\ \forall j<i\;\; w(i)<w(j) \}.
\]
\end{df}
 For example let $w=[3,2,7,8,4,6,1,5]$, then $\{2,7\}$ are 
the l.t.r.min according
to~\ref{df4}.1, and $\{1,2\}$ according to~\ref{df4}.2.\\
\newline
With either definition we have
%
\begin{pro}\label{pro4'}
Let $w\in S_n$, 
then $del_S(w)$ equals the number of l.t.r.min of $w^{-1}$
according to either Definition~\ref{df4}.1 or~\ref{df4}.2. Since 
by Lemma~\ref{lem9} $del_S(w)=del_S(w^{-1})$, this also equals 
the number of l.t.r.min of $w$.
In particular, 
\[
\mid Del_S(w)\mid =del_S(w)=del_S(w^{-1}).
\]
\end{pro}
%
%
%
\noindent
{\it  Proof.} By induction on $n\ge 2$. 
First, $S_2=\{1,s_1\}$ and $s_1=[2,1]$ has one l.t.r.min - according to 
either~\ref{df4}.1 or~\ref{df4}.2.  
Proceed now with the inductive step, which is essentially the same for
both definitions.
Let $w=w_1\cdots w_{n-1}$ be 
the canonical presentation of $w$, let $\sg =  w_1\cdots w_{n-2}$ 
(so $\sg\in S_{n-1}$)
and assume true for $\sg$. Write $\sg^{-1}=[b_1,\ldots , b_{n-1},n]$.
If $ w_{n-1}=1$, the proof is given by the induction hypothesis. Otherwise 
$ w^{-1}_{n-1}= s_ks_{k+1}\cdots s_{n-1}$ for some 
$1\le k\le n-1$. Denoting 
$s_{[k,n-1]}= s_ks_{k+1}\cdots s_{n-1}$ we see that
$w^{-1}=s_{[k,n-1]}\sg ^{-1}$. 
Comparing $\sg ^{-1}$ with $w^{-1}=s_{[k,n-1]}\sg ^{-1}$, we see that 
\begin{enumerate}
\item
the (position with) $n$ in $\sg ^{-1}$ is replaced in  $w^{-1}$  by $k$;
\item
each $j$ in $\sg ^{-1}$, $k\le j\le n-1$, is replaced by $j+1$ in $w^{-1}$;
\item
each $j$, $1\le j \le k-1$ is unchanged. 
\end{enumerate}

Thus $\sg^{-1}=[b_1,\ldots , b_{n-1},n]$, 
$w^{-1}=[c_1,\ldots , c_{n-1},k]$, and the two tuples
$(b_1,\ldots , b_{n-1})$ and $(c_1,\ldots , c_{n-1})$ are order-isomorphic.
This implies that if $k > 1$ then $\sg ^{-1}$ and $w^{-1}$ have the 
same left-to-right minima. Let $k=1$ and adopt Definition~\ref{df4}.1 first, 
then $w^{-1}$ has $i=n$ as
an additional left-to-right minimum, which  
completes the proof in that case.  In the
case of Definition~\ref{df4}.2, compare $k=2$ with $k=1$ to 
deduce that $i$ such
that $c_i=2$ is an additional l.t.r.min, and the proof follows.
\qed
\begin{rem}\label{rem5}
The above proof implies a bit more: 
Note that the above case $k=1$ is equivalent to both $n\in Del_S(w^{-1})$
and to $\ep_{S,n-1}(w)=1$, where $\ep_{S,i}(w)$ 
are given by  Definition~\ref{df5}. By induction on $n$, the above proof
implies that $Del_S(w^{-1})=\{i+1\mid \ep_{S,i}(w)=1\}$.
Let now $D\subseteq [n-1]$ and let $\pi\in S_n$. 
The condition $D=Del_S(\pi ^{-1})$ implies that
$D=\{i+1\mid\ep _{S,i}(\pi)=1\}$; this uniquely determines 
${\bar\ep}_S(\pi)$, and hence determines a unique value 
$t^{\ep _D}:= t^{{\bar\ep}_S(\pi)}$: if $D\ne H$ then 
$t^{\ep _D}\ne t^{\ep _H}$. We shall apply this observation in the
proof of Theorem~\ref{shuff.1}.
\vskip 0.1 truecm
 
\vskip 0.1 truecm
%
%
%
%
\end{rem}

Each of the two definitions of l.t.r.min can be extended as follows.
\begin{df}\label{df6}
Let $w=[b_1,\ldots ,b_{n}]\in S_n$. Then $1\le i \le n$ is an
a.l.t.r.min (almost-left-to-right minimum) if, first of all, there is at most 
one $b_j$ smaller than $b_i$ and left of $b_i$:
card$\{1\le j \le i \mid b_j < b_i \} \le 1$.\\
The second condition is one of the following:\\
Either\\
1. Do not count $i=1$ and $i=2$  as a.l.t.r.min
(which is Definition~\ref{df6}.1 of a.l.t.r.min), \\
or:\\
2. Do not count $i$ such that $b_i=1,2$ as an a.l.t.r.min
(which is Definition~\ref{df6}.2 of a.l.t.r.min).\\
3. For $w\in A_{n+1}$ define $Del_A(w)$ to be the set of a.l.t.r.min
of $w$ according to  Definition~\ref{df6}.1. 

\end{df}
\begin{rem}\label{rem3}
\begin{enumerate}
\item
Without the above restrictions 1 and 2 in Definition~\ref{df6},
$i$ such that $b_i\in \{1,2\}$
is an  a.l.t.r.min; 
similarly, if $i\in\{1,2\}$ then $i$ is an a.l.t.r.min. 
\item
If $b_i=1$ and $b_j=2$ are interchanged in $w=[b_1,\ldots ,b_{n}]$, this does
not change the set of a.l.t.r.min indices. Also, if $b_1$ and $b_2$ 
are interchanged this would not change the set of a.l.t.r.min indices.
Thus with either definition~\ref{df6}.1 or ~\ref{df6}.2,
$s_1w$ and $ws_1$ have the same a.l.t.r.min as $w$ itself.
\end{enumerate}
\end{rem}
%
%
\begin{pro}\label{pro9}
Let $w\in S_n$, then the number of occurrences of $s_2$ in (the canonical presentation of)
$w$ equals the number of a.l.t.r.min of $w^{-1}$ according to either 
Definition~\ref{df6}.1 or~\ref{df6}.2. Lemma~\ref{lem9} 
implies that this is also the number of a.l.t.r.min of $w$.
\end{pro}
{\it Proof.} By induction on $n$. This is easily verified for 
$n=2$, and we proceed with the inductive step. \\
Let $w=w_1\cdots w_{n-1}$ be the canonical presentation of $w$, and denote 
$\sg =w_1\cdots w_{n-2}$, so that $w^{-1}=w^{-1}_{n-1}\sg^{-1}$.
If $w_{n-1}=1$ we are done by induction. Otherwise, by the S-procedure, 
$w_{n-1}=s_{n-1}\cdots s_kx$ where $k\ge 2$ and $x\in\{1,s_1\}$.\\
Write $w^{-1}=xs_k\cdots s_{n-1}\sg ^{-1}=xs_{[k,n-1]}\sg ^{-1}$.
By Remark~\ref{rem3}, 
$s_{[k,n-1]}\sg ^{-1}$ and $xs_{[k,n-1]}\sg ^{-1}$ 
have the same number of a.l.t.r.min.
Therefore it suffices to show:
\begin{enumerate}                                                         
\item
If $k\ge 3$ then $\sg ^{-1}$  has equal number
of a.l.t.r.min as $s_{[k,n-1]}\sg ^{-1}$.
\item 
If $k=2$, $s_{[2,n-1]}\sg ^{-1}$ has one more a.l.t.r.min than
$\sg ^{-1}$.
\end{enumerate}
Let $\sg ^{-1}=[b_1,\ldots ,b_{n-1},n]$, then 
$s_{[k,n-1]}\sg ^{-1}=[c_1,\ldots ,c_{n-1},k]$, and as in the proof 
of~\ref{pro4'}, $(b_1,\ldots ,b_{n-1})$ 
and $(c_1,\ldots ,c_{n-1})$ are order isomorphic . 
If $k\ge 3$, the last position (with $k$) is not an a.l.t.r.min, while if
$k=2$, it is an additional a.l.t.r.min, and this implies the proof
in the case of~\ref{df6}.1.  In case of~\ref{df6}.2, compare $k=3$
with $k=2$: a $2$ is changed into a $3$, which is an additional  
a.l.t.r.min.  
\qed\\
\newline
By essentially the same argument we have
\begin{pro}\label{pro8}
Let $v\in A_{n+1}$ then, with either 
Definition~\ref{df6}.1 or~\ref{df6}.2 
of a.l.t.r.min, 
$del_A(v)$ equals the number of a.l.t.r.min of $v^{-1}$. 
In particular $\mid Del_A(v)\mid = del_A(v)= del_A(v^{-1})$.
\end{pro}
{\it Proof.} Again, by induction on $n$. This is easily verified for 
$n+1=3$, so proceed with the inductive step. \\
Let $v=v_1\cdots v_{n-1}$ be the A- canonical presentation of $v$, and denote 
$\sg =v_1\cdots v_{n-2}$, so that $v^{-1}=v^{-1}_{n-1}\sg^{-1}$.
If $v_{n-1}=1$ we are done by induction. Otherwise, by the A-procedure, 
$v_{n-1}=xs_n\cdots s_ky$ where $k\ge 2$ and $x,y\in\{1,s_1\}$;
moreover, $k=2$ if and only if either $a_1$ or $a^{-1}_1$ occurs in 
$v_{n-1}$.\\
Write $v^{-1}=ys_k\cdots s_nx\sg ^{-1}=ys_{[k,n]}x\sg ^{-1}$, and proceed
as in the proof of~\ref{pro9}, applying~\ref{rem3}.2. \qed
%
%
%

%
\begin{rem}\label{rem4}
Given $w\in S_n$,
one can define a.a.l.t.r.min, a.a.a.l.t.r.min,  etc, then prove the 
corresponding propositions, analogue of Proposition~\ref{pro9}. For 
example, we have
\begin{df}\label{df7}
Let $w=[b_1,\ldots ,b_{n}]\in S_n$. Then $1\le i \le n$ is an
a.a.l.t.r.min (almost-almost-left-to-right minimum) if
card$\{1\le j \le i \mid b_j < b_i \} \le 2$ and\\ 
1. $i\ne 1,\,2,\,3$ (which is Definition~\ref{df7}.1
of a.a.l.t.r.min),\\
or\\
2. $b_i\ne 1,\,2,\,3$.  (which is Definition~\ref{df7}.2 of a.a.l.t.r.min).
%
\end{df}
One can then prove that, with either definition of a.a.l.t.r.min,
the number of a.a.l.t.r.min of $w\in S_n$ equals
the number of occurrences of $s_3$ in $w$. Similarly for the occurrences 
of the other $s_i$'s.
\end{rem}
Similar to Proposition~\ref{pro3}, we define $w_S(n,\ell,k)$ to be the 
number of $S$ canonical words in $S_n$ with $\ell$ occurrences of $s_k$
(define $w_A(n+1,\ell,k)$ similarly), and we have
\begin{pro}\label{pro10}
Let $k\leq n-1$, then 
\[
\sum _{\ell=0}^{n-k}w_S(n,\ell,k)t^\ell=k!(kt+1)(kt+2)\cdots (kt+n-k),
\]
hence $w_S(n,\ell,k)=k!k^{\ell}c(n-k+1,\ell +1)$, 
and similarly for $w_A(n+1,\ell,k)$.
\end{pro}
{\it Proof} -  is omitted.

\section{Lemmas on Shuffles}

In this section we prove lemmas, which will be used in the next section
to prove the main theorem.

\subsection{Equi-distribution on Shuffles}

The following result follows from Theorem~\ref{p3}.

\begin{pro}\label{1}
Let $i\in [n-1]$, and let
$\pi\in S_n$ with 
$\supp(\pi)\subseteq [i]$.
Then
$$
\sum\limits_{\Des_S(r^{-1})\subseteq \{i\}}
q^{\rmaj_{S_n}(\pi r)-\rmaj_{S_i}(\pi)}=
\sum\limits_{\D(r^{-1})\subseteq \{i\}} 
q^{\ell_S(\pi r)-\ell_S(\pi)}=
{n\brack i}_q. 
$$
\end{pro}
\noindent{\bf Proof.}
Let $\rho_n:=(1,n)(2,n-1)\cdots\in S_n$
and $\rho_i:=(1,i)(2,i-1)\cdots\in S_i$.
By (\ref{eq0})
$$
\sum\limits_{\Des_S(r^{-1})\subseteq \{i\}}
q^{\rmaj_{S_n}(\pi r)-\rmaj_{S_i}(\pi)}=
\sum\limits_{\Des_S(r^{-1})\subseteq \{i\}}
q^{maj_{S}(\rho_n\pi r\rho_n)-maj_{S}(\rho_i\pi\rho_i)}=
$$
$$
\sum\limits_{\Des_S(r^{-1})\subseteq \{i\}}
q^{maj_{S}(\rho_n\pi \rho_n \rho_n r\rho_n)-maj_{S}(\rho_i\pi\rho_i)}=
\sum\limits_{\Des_S(\hat r^{-1})\subseteq \{n-i\}}
q^{maj_{S}(\rho_n\pi \rho_n \hat r)-maj_{S}(\rho_i\pi\rho_i)}.
$$
The last equality follows from (\ref{eq011}).

\noindent
Note that $\supp(\rho_n \pi\rho_n)\subseteq [n-i+1,n]$
and verify that
$\nu_{n-i}^{-1} \rho_n \pi\rho_n\nu_{n-i}=
\rho_i \pi\rho_i$, where $\nu_{n-i}:=(1,n-i+1)(2,n-i+2)\cdots$.
Indeed, let $j\le i$, then $\nu_{n-i}(j)=j+n-i$, hence
$\rho_n\nu_{n-i}(j)=\rho_n (j+n-i)=n-(j+n-i)+1=i-j+1=\rho_i (j)$. Similarly,
if $k\le i$, also $\nu^{-1}_{n-i}\rho_n (k) = \rho_i (k)$. This implies the
above equality.
Now, obviously $\supp(1)\subseteq [n-i]$ and $maj_S(1)=0$.
Thus by Garsia-Gessel's Theorem (Theorem
\ref{p3}) (taking $\pi_1=1$ and $\pi_2=\rho_n\pi\rho_n$)
the right-hand-side is equal to
$$
\sum\limits_{\Des_S(\hat r^{-1})\subseteq \{n-i\}}
q^{maj_{S}(1\cdot\rho_n\pi \rho_n \cdot\hat r)-maj_S(1)-
maj_{S}(\nu_{n-i}^{-1}\rho_n\pi\rho_n\nu_{n-i})}={n\brack i}_q. 
$$

\smallskip

The equality
$$
\sum\limits_{\D(r^{-1})\subseteq \{i\}} 
q^{\ell_S(\pi r)-\ell_S(\pi)}=
\sum\limits_{\D(r^{-1})\subseteq \{i\}} 
q^{\ell_S(\pi\cdot 1\cdot r)-\ell_S(\pi)-\ell_S(1)}
={n\brack i}_q 
$$
is an immediate consequence of Fact \ref{p3i}, combined with (\ref{eqn3i}).

\qed

\medskip

\begin{nt}\label{nt1}
Let $r$ be an $\{i\}$-shuffle and let $supp(\pi)\subseteq [i]$ as above.
If $r(1)\ne 1$, necessarily $r(1)=i+1$,
hence also $\pi r (1)=i+1$. It follows that 
$$
\pi r (1)\in \{\pi(1), i+1\}.
$$
\end{nt}

The next lemma requires some preparations. 

Fix $1\le i\le n-1$ and define $g_i:S_n\to S_{n-1}$ as follows: Let \\
$\sg=[a_1,\ldots ,a_n]\in S_n$, then $g_i(\sg)=[a'_1,\ldots ,a'_{n-1}]$
is defined as follows: delete $a_j=i+1$, leave $a'_k=a_k$ unchanged if
$a_k\le i$, and change $a'_t=a_t -1$ if $a_t\ge i+2$. Denote
 $g_i(\sg)=\sg '$. For example, let $\sg =[5, 2,3,6,1,4]$ and $i=2$,
then $g_2(\sg)=\sg '=[4,2,5,1,3]$.
Let $supp(\pi )\subseteq [i]$, then $g_i(\pi)=\pi$: $\pi '=\pi$.
Moreover, since $\pi$ only permutes $1,\ldots , i$, the following 
basic property of $g_i$ is rather obvious, 
since $supp(\pi)\subseteq [i]$:

\begin{fac}\label{gi}
1. Let $\sg\in S_n$, then $\pi (g_i \sg )=g_i(\pi\sg )$, namely,
$(\pi\sg )'=\pi ' \sg' =\pi\sg '$. \\

\smallskip

2. $g_i$ is a bijection between the $\{i\}$-shuffles $r\in S_n$ 
satisfying $r(1)=i+1$,
and all the $\{i\}$-shuffles $r'\in S_{n-1}$:
\[
g_i:\{r\in S_n\mid Des_S(r^{-1})\subseteq \{i\},\; r(1)=i+1\}
\to \{r'\in S_{n-1}\mid Des_S(r^{-1})\subseteq \{i\}\}
\]
is a bijection.
\end{fac}
\vskip 0.15 truecm
\noindent We need
\begin{lem}\label{g}
Let $1\le i\le n-2$, $supp(\pi)\subseteq [i]$ and $r(1)=i+1$. Also
let $g_i (\pi)=\pi '$ and $g_i (r)=r '$.
\begin{enumerate}
\item
If $r(2)=i+2$ then $\rmaj_{S_n}(\pi r)=\rmaj_{S_{n-1}}(\pi' r')$.
\item
If $r(2)=1$ then $\rmaj_{S_n}(\pi r)=n-1+\rmaj_{S_{n-1}}(\pi' r')$.
\end{enumerate}
\end{lem}
{\it Proof.} 
By Note \ref{nt1}, 
$\pi r=[i+1,a_2,\ldots , a_n]$ then, applying $g_i$, we have 
$\pi ' r'=[a'_2,\ldots , a'_n]$, and it is easy to check that for all
$2\le k\le n-1$, ~~$a_k > a_{k+1}$ if and only if $a'_k > a'_{k+1}$.
Thus, for $2\le k\le n-1$, ~~$k\in Des(\pi r)$ ~~if 
and only if ~~$k-1\in Des(\pi' r')$; 
note also that such $k$ contributes $n-k=(n-1)-(k-1)$ 
to both $\rmaj_{S_n}(\pi r)$
and to $\rmaj_{S_{n-1}}(\pi' r')$.

\begin{enumerate}
\item
If $r(2)=i+2$ then $a_2=\pi r  (2)=i+2$, hence $1\not\in Des(\pi r)$,
and the descents of $\pi r$ occur only for (some) $2\le k\le n-1$, 
and the above argument implies the proof.
\item
If $r(2)=1$ then $a_2=\pi r (2) =\pi (1)< i+1$, hence $1$ is a descent
of $\pi r$, contributing $n-1$ to $\rmaj_{S_n}(\pi r)$, and again, 
the above argument completes the proof.\qed
\end{enumerate} 
%

\begin{lem}\label{22}
With the notations of Proposition \ref{1} 
$$
\sum\limits_{\Des_S(r^{-1})\subseteq \{i\}\hbox{ and } \pi r(1)=i+1}
q^{\rmaj_{S_n}(\pi r)-\rmaj_{S_i}(\pi)}=q^i{n-1\brack i}_q.  \leqno(1)
$$
and
$$
\sum\limits_{\Des_S(r^{-1})\subseteq \{i\}\hbox{ and } \pi r(1)=\pi(1)}
q^{\rmaj_{S_n}(\pi r)-\rmaj_{S_i}(\pi)}={n-1\brack i-1}_q.  \leqno(2)
$$
\end{lem}

\noindent{\bf Proof.}
By induction on $n-i$. For $n-i=1$, the $\{n-1\}$-shuffles are
$[1,\ldots,j-1,n,j,\ldots,n-1]=[1,\ldots,n]s_{n-1}s_{n-2}\cdots s_j$,
$\;1\le j\le n-1$.
Thus the summation in (2) is over
$r\in R_{n-1}^S - \{s_{n-1} s_{n-2}\cdots s_1\}$ and 
Equation (2) follows from Remark~\ref{rem2} (with $n-1$ replacing $n$).
Now, 
\[
sum (1)+sum (2)=\sum\limits_{\Des_S(r^{-1})\subseteq \{n-1\}}
q^{\rmaj_{S_n}(\pi r)-\rmaj_{S_{n-1}}(\pi)},
\] 
Hence, by Proposition~\ref{1} 
$$ sum (1)+sum (2)={n\brack n-1}_q,$$
so 
$$sum(1)={n\brack n-1}_q -{n-1\brack n-2}_q=q^{n-1},$$
which verifies (1) in that case.
\vskip 0.1 truecm

Let now  $n-i\ge 2$ and assume the lemma holds for $n-1-i$. 
\begin{description}
\item[(1)]
Since $Des(r^{-1})\subseteq \{i\}$ and $r(1)=i+1$,
either $r(2)=i+2$ (then $\pi r(2)=i+2$), or $r(2)=1$ (then $\pi r(2)=\pi (1)$).
%
Thus, the sum in (1) equals $sum[r(2)=i+2]+sum[r(2)=1]$. Apply $g_i$ 
to the permutations in these sums, and apply Lemma~\ref{g}.1
and Fact~\ref{gi}; then, by induction on $n$,
$$
sum[r(2)=i+2]=
\sum\limits_{\Des_S(r'^{-1})\subseteq \{i\}\hbox{ and } \pi' r'(1)=i+1}
q^{\rmaj_{S_{n-1}}(\pi' r')-\rmaj_{S_i}(\pi')}
$$
$$
=q^i{n-2\brack i}_q. 
$$
Similarly, by Lemma~\ref{g}.2 and Fact~\ref{gi}, 

$$
sum[r(2)=1]=
\sum\limits_{\Des_S(r'^{-1})\subseteq \{i\}\hbox{ and } \pi' r'(1)=\pi'(1)}
q^{n-1+\rmaj_{S_{n-1}}(\pi' r')-\rmaj_{S_i}(\pi')}
$$
$$
=q^{n-1}{n-2\brack i-1}_q.  
$$
%
%
%
Adding the last two sums, we conclude:
$$
\sum\limits_{\Des_S(r^{-1})\subseteq \{i\}\hbox{ and } \pi r(1)=i+1}
q^{\rmaj_{S_{n-1}}(\pi r)-\rmaj_{S_i}(\pi)}=q^i{n-2\brack i}_q+q^{n-1}{n-2\brack i-1}_q=
$$
$$
=q^i({n-2\brack i}_q+q^{n-1-i}{n-2\brack i-1}_q)=
q^i{n-1\brack i}_q.
$$
\item[(2)] is an immediate consequence of Proposition \ref{1} and part (1), since
$$
{n\brack i}_q-{n-1\brack i-1}_q=q^i{n-1\brack i}_q.
$$
\end{description}

\qed

We have an analogous lemma for length.

\begin{lem}\label{2}
With the notation of Proposition \ref{1}
$$
\sum\limits_{\Des_S(r^{-1})\subseteq \{i\}\hbox{ and } \pi r(1)=i+1}
q^{\ell_S(\pi r)-\ell_S(\pi)}=
q^i{n-1\brack i}_q. \leqno(1)
$$
and
$$
\sum\limits_{\Des_S(r^{-1})\subseteq \{i\}\hbox{ and } \pi r(1)=\pi(1)}
q^{\ell_S(\pi r)-\ell_S(\pi)}=
{n-1\brack i-1}_q. \leqno(2)
$$
\end{lem}

\noindent{\bf Proof}. The case $n-i=0$ is obvious  (the sum in (1)
is empty while in (2), $r = 1$), so assume 
$i\le n-1$. 
Recall that in general, $\ell _S(\sg)$ equals the
number $inv_S(\sg)$ of inversions of $\sg$. 

We prove (1) first, so let $\pi r (1)=i+1$. As in Lemma~\ref{g},
write 
\[
\pi r =[i+1,a_2,\ldots ,a_n]\quad\mbox{and}\quad
g_i(\pi r)=\pi' r'=[a'_2,\ldots ,a'_n],
\]
and compare their inversions. Clearly, $i+1$ contributes $i$ inversions
to $inv_S(\pi r)$. 
Also, as in the proof of Lemma~\ref{g}, 
there is a bijection between the inversions among $\{a_2,\ldots , a_n\}$
and those among  $\{a'_2,\ldots , a'_n\}$. Thus 
$inv_S(\pi r)=i+inv_S(\pi' r')$. 
Also, since $supp(\pi)\subseteq [i]$ , $inv_S(\pi )=inv_S(\pi' )$.
Induction, Fact~\ref{gi}
and Proposition~\ref{1} imply the proof of (1).  
Now, by Proposition~\ref{1}, (1) implies the proof of (2). \qed

\subsection{Canonical Presentation of Shuffles}

\begin{obs}\label{der6'}
Let $1\le i<n$. Every $\{i\}$-shuffle
has a unique canonical presentation of the form
$w_i w_{i+1}\cdots w_{n-1}$, where $\ell(w_j)\ge \ell(w_{j+1})$
for all $j\ge i$.
\end{obs}
{\it Proof.} Apply the `S-Procedure' that follows Theorem~\ref{thm1}. 
Note that 
after pulling $n,n-1,\ldots , i+1$ to the right, an $\{i\}$-shuffle
is transformed into the identity permutation.
\qed\\
\newline
Let $\bar\ep=(\ep _1, \ldots , \ep _{n-1})$, then denote
$t^{\bar\ep}=t^{\ep _1}_1\cdots t^{\ep _{n-1}}_{n-1}$.
\begin{cor}\label{der6''} Recall Definition~\ref{df5}.
For an $\{i\}$-shuffle $w$,
$$
\del_S(w)=\cases
{1, & \hbox{ if\ $w(1)=i+1$;}\cr
0, & \hbox { otherwise, }\cr}
$$
and therefore
$$
t^{{\bar\ep}_S(w)}=t^{del_S(w)}_i=\cases
{t_i, & \hbox{ if\ $w(1)=i+1$;}\cr
1, & \hbox { otherwise }\cr}
$$
\end{cor}
\noindent{\it Proof.} Write $w=w_iw_{i+1}\cdots w_{n-1}$ (canonical presentation) with 
$\ell _S(w_i)\ge\cdots\ge\ell _S(w_{n-1})$, then $\ep_{S,j}(w)=0$ for $j>i$.
Thus $del_S(w)$ is either 1 or 0, and is 1 exactly when 
$w_i=s_i\cdots s_1$, in which case $w(1)=i+1$.
\qed
\begin{rem}\label{rem6}
Let $r,\pi\in S_n$, $r$ an $\{i\}$-shuffle and $supp(\pi)\subseteq [i]$.
Then the corresponding
canonical presentations are: $\pi = w_1\cdots w_i,\quad$ 
$r= w_{i+1}\cdots w_{n-1}$, hence also $\pi r= w_1\cdots w_{n-1}$ is canonical presentation.
In particular, 
${\bar\epsilon_S(\pi r)}={\bar\epsilon_S(\pi)+{\bar\epsilon_S(r)}}$.

\end{rem}

\vskip 0.3 truecm
We generalize:
Let $B=\{i_1,i_2\}$ and let $w\in S_n$ be a $B$-shuffle. Then $w$ shuffles  
the three subsets $\{1,\ldots ,i_1\}$, $\{i_1+1,\ldots ,i_2\}$ and
$\{i_2+1,\ldots ,n\}$. Clearly $w$ 
has a unique presentation as a 
product $w=\tau_1\tau_2$ where
$\tau_2\in S_n$ shuffles $\{1,\ldots ,i_2\}$ with $\{i_2+1,\ldots ,n\}$, and 
$\tau_1\in S_{i_2}$ shuffles $\{1,\ldots ,i_1\}$ with $\{i_1+1,\ldots ,i_2\}$.
By Observation~\ref{der6'}, $\tau_1 =w_{i_1}w_{i_1+1}\cdots w_{i_2-1}$ and
$\tau_2 =w_{i_2}w_{i_2+1}\cdots w_{n-1}$, where each $w_j\in R^S_j$.
Thus 
\[
w=w_{i_1}\cdots w_{i_2-1}w_{i_2}\cdots w_{n-1}
\]
is the $S$ canonical presentation of $w$, 
\[
del_S(w)=del_S(\tau_1) + del_S(\tau_2)\qquad\mbox{and}\qquad
t^{{\bar\ep}_S(w)}=t^{del_S(\tau _1)}_{i_1}t^{del_S(\tau _2)}_{i_2}.
\]
This easily generalizes to an arbitrary 
$B=\{i_1,\dots,i_k\}\subseteq \{1,\dots,n-1\}$, which proves the 
following proposition.
\begin{pro}\label{pro14}
Let $B=\{i_1,\dots,i_k\}\subseteq \{1,\dots,n-1\}$ and let $i_{k+1}:=n$.
Every $B$-shuffle  $\pi\in S_n$ has a unique presentation
$$
\pi=\tau _1\cdots \tau _k
$$
where $\tau_j$ is an $\{i_j\}$-shuffle in $S_{i_{j+1}}$ (for $1\le j\le k$).
Moreover,
\[
\del_S(\pi)=\sum\limits_{j=1}^k \del_S(\tau _j)\qquad\mbox{and}\qquad
t^{{\bar\ep}_S(\pi)}=t^{del_S(\tau _1)}_{i_1}\cdots
t^{del_S(\tau _k)}_{i_k}.
\]
\end{pro}
\vskip 0.15 truecm
%
%
%
%
%
\medskip
%

%
%
\section{The Main Theorem}

\vskip 0.1 truecm
Recall the definitions of the $A$-descent set
$\Des_{A}$ and the $A$-descent number $\des_A$
(Definition~\ref{df8}).
Let $B\subseteq [n-1]$ and $\pi\in S_n$. Recall from Fact~\ref{p3.1}
that $Des_S(\pi ^{-1})\subseteq B$ if and only if $\pi$ is a $B$-shuffle.

\vskip 0.1 truecm

\vskip 0.1 truecm
The following is our main theorem, which we now prove.


%
%
%
%
%
%
\begin{thm}\label{shuff.1}
For every subsets $D_1\subseteq [n-1]$ and $D_2\subseteq [n-1]$
$$
\sum\limits_{\{\pi\in S_n|\ \Des_S(\pi^{-1})\subseteq D_1,\ 
\Del_S(\pi^{-1})\subseteq D_2\}} q^{\rmaj_{S_n}(\pi)}
=\leqno(1)
$$
$$
\sum\limits_{\{\pi\in S_n|\ \Des_S(\pi^{-1})\subseteq D_1,\
\Del_S(\pi^{-1})\subseteq D_2\}} q^{\ell_S(\pi)},
$$
and
$$
\sum\limits_{\{\sigma\in A_{n+1}|\  \Des_A(\sigma^{-1})\subseteq D_1,\
\Del_A(\sigma^{-1})\subseteq D_2\}} q^{\rmaj_{A_{n+1}}(\sigma)}
=\leqno(2)
$$
$$
\sum\limits_{\{\sigma\in A_{n+1}|\  \Des_A(\sigma^{-1})\subseteq D_1,\
\Del_A(\sigma^{-1})\subseteq D_2\}} q^{\ell_A(\sigma)}.
$$
\end{thm}

\vskip 0.1 truecm

An immediate consequence of Theorem~\ref{shuff.1} is

\begin{cor}\label{shuff.2}
$$
\sum\limits_{\pi\in S_n} q_1^{\rmaj_{S_n}(\pi)}
q_2^{\des_S(\pi^{-1})}q_3^{\del_S(\pi^{-1})}=
\sum\limits_{\pi\in S_n} q_1^{\ell_S(\pi)}
q_2^{\des_S(\pi^{-1})}q_3^{\del_S(\pi^{-1})}.
\leqno(1)
$$
$$
\sum\limits_{\sigma\in A_n} q_1^{\rmaj_{A_{n+1}}(\sigma)} q_2^{\des_A(\sigma^{-1})}q_3^{\del_A(\sigma^{-1})}=
\sum\limits_{\sigma\in A_n} q_1^{\ell_A(\sigma)}
q_2^{\des_A(\sigma^{-1})}q_3^{\del_A(\sigma^{-1})}.
\leqno(2)
$$
\end{cor}


\medskip

Note that in Corollary \ref{shuff.2}(1) both definitions \ref{df4}.1
and \ref{df4}.2 for calculating
 $\del_S$ could be used. This follows from Proposition~\ref{pro4'}.
Similarly, in Corollary \ref{shuff.2}(2)
both definitions \ref{df6}.1
and \ref{df6}.2 for calculating $\del_A$ could be used 
(by Proposition~\ref{pro8}).

\subsection{A Lemma}

\begin{lem}\label{der7} 
Let $i\in [n]$, and let
$\sg$ be a permutation in $S_n$, such that 
$\supp(\sg)\subseteq [i]$.
Then
$$
\sum\limits_{\D(r^{-1})\subseteq \{i\}} 
q^{\ell_S(\sg r)}t^{\bar\epsilon_S(\sg r)}=
q^{\ell_S(\sg)}t^{\bar\epsilon_S(\sg)}\cdot
\left({n-1\brack i-1}_q+t_i q^{i}{n-1\brack i}_q\right) \leqno (1)
$$
and
$$
\sum\limits_{\D(r^{-1})\subseteq \{i\}}
q^{\rmaj_{S_n}(\sg r)}t^{\bar\epsilon_S(\sg r)}=
q^{\rmaj_{S_i}(\sg)}t^{\bar\epsilon_S(\sg)}\cdot
\left({n-1\brack i-1}_q+t_i q^{i}{n-1\brack i}_q\right). \leqno (2)
$$
\end{lem}

\noindent{\bf Proof.} 
By Definition~\ref{df5} and Remark~\ref{rem6}

$$
t^{\bar\epsilon_S(\sg r)}=t^{\bar\epsilon_S(\sg)+{\bar\epsilon_S(r)}},
$$
and by Corollary~\ref{der6''} 
$$
t^{{\bar\ep}_S(r)}=\cases
{t_i, & \hbox{ if\ $r(1)=i+1$;}\cr
1, & \hbox { otherwise }\cr}
$$
Noting that $r(1)=i+1$ if and only if $ \sg r(1)=i+1$, 
and recalling that $\sg r(1)\in\{\sg (1),i+1\}$, we obtain
$$
t^{{\bar\ep}_S(\sg r)}=\cases
{t^{{\bar\ep}_S(\sg)}t_i, 
& \hbox{ if\ $\sg r(1)=i+1$;}\cr
t^{{\bar\ep}_S(\sg)}, 
& \hbox{ if\ $\sg r(1)=\sg(1)$}\cr}
$$
Combining this with Lemmas \ref{22} and \ref{2} gives the desired result. 
For example, regarding length,
$$
\sum\limits_{\D(r^{-1})\subseteq \{i\}} 
q^{\ell_S(\sg r)}t^{\bar\epsilon_S(\sg r)}=
$$
$$
=\sum\limits_{\D(r^{-1})\subseteq \{i\}\hbox{ and } \sg r(1)=\sg(1)} 
q^{\ell_S(\sg r)}t^{\bar\epsilon_S(\sg r)}+
\sum\limits_{\D(r^{-1})\subseteq \{i\}\hbox{ and } \sg r(1)=i+1} 
q^{\ell_S(\sg r)}t^{\bar\epsilon_S(\sg r)}=
$$ 
$$
=q^{\ell_S(\sg)}t^{\bar\epsilon_S(\sg)}\cdot
\left({n-1\brack i-1}_q+t_i q^{i}{n-1\brack i}_q\right).
$$
This proves part (1). A similar argument proves (2).

\qed

\subsection{Proof of Main Theorem}
%
%
%
%
%



\bigskip

\noindent{\bf Proof of Theorem \ref{shuff.1}(1).}

\noindent
By the principle of inclusion and exclusion, we may replace 
$\Del_S(\pi^{-1})\subseteq D_2$  by $\Del_S(\pi^{-1})= D_2$
in both 
hand-sides of Theorem \ref{shuff.1}(1).
By Remark~\ref{rem5}, 
$\{\pi\in S_n\ |\ \Del_S(\pi^{-1})= D_2 \}$ (i.e. the set $D_2$)
%
%
determines the 
unique value  $t^{\ep_{D_2}}:=   t^{\bar\ep_S(\pi)} $.
\vskip 0.1 truecm

\noindent Hence, Theorem \ref{shuff.1}(1) is equivalent to the following statement :

For every subset $B\subseteq [n-1]$ 
$$
\sum\limits_{\{\pi\in S_n|\ \Des_S(\pi^{-1})\subseteq B\}} 
q^{\rmaj_{S_n}(\pi)}t^{\bar\ep_S(\pi)}
=
$$
$$
\sum\limits_{\{\pi\in S_n|\ \Des_S(\pi^{-1})\subseteq B\}} q^{\ell_S(\pi)}t^{\bar\ep_S(\pi)},
$$

\medskip

This statement is proved by induction on the cardinality of $B$.
If $|B|=1$ then $B=\{i\}$ for some $i\in [n-1]$ and 
Theorem \ref{shuff.1}(1)  is given by 
Lemma \ref{der7} (with $\sg=1$).
Assume that the theorem holds for every $B\subseteq [n-1]$ of cardinality less than $k$.
Let $B=\{i_1,\dots,i_k\}\subseteq [n-1]$
and denote $\bar B:=\{i_1,\dots,i_{k-1}\}$. By Proposition \ref{pro14},
for every $\pi\in S_n$ with $\Des_S(\pi^{-1})\subseteq B$
there is a unique presentation
$$
\pi=\bar\pi \tau _k,
$$
where $\bar\pi$ is a $\bar B$-shuffle in $S_{i_k}$ and
$\tau_k$ is an $\{i_k\}$-shuffle in $S_{n}$.
Moreover, $\Des_S(\pi^{-1})\subseteq B$
if and only if $\pi$ has such a presentation. Hence
$$
\sum\limits_{\{\pi\in S_n|\ \Des_S(\pi^{-1})\subseteq B\}} 
q^{\rmaj_{S_n}(\pi)}t^{\bar\ep_S(\pi)}=
$$
$$
\sum\limits_{\{\bar\pi\in S_{i_k}, \tau_k\in S_n|\ \D_S(\bar\pi^{-1})
\subseteq \bar B, 
\D_S(\tau_k^{-1})\subseteq\{i_k\}\}} 
q^{\rmaj_{S_n}(\bar\pi\tau_k)}t^{\bar
\ep_S(\bar\pi\tau_k)}=
$$
$$
\sum\limits_{\{\bar\pi\in S_{i_k}|\ \D_S(\bar\pi)\subseteq\bar B\}}\ \ \ 
\sum\limits_{\{\tau_k\in S_n|\ \D_S(\tau_k^{-1})\subseteq \{i_k\}\}}
q^{\rmaj_{S_n}(\bar\pi\tau_k)}t^{\bar
\ep_S(\bar\pi\tau_k)}.
$$
By Lemma \ref{der7}(2) this equals to\\

\vskip 0.1 truecm

$$
\sum\limits_{\{\bar\pi\in S_{i_k}|\ \D_S(\bar\pi^{-1})\subseteq \bar B\}}
q^{\rmaj_{S_{i_{k-1}}}(\bar\pi)}t^{\bar\epsilon_S(\bar\pi)}\cdot
\left({n-1\brack i_k-1}_q+t_{i_k} q^{i}{n-1\brack i_k}_q\right)
$$
which, by induction, equals
$$
\sum\limits_{\{\bar\pi\in S_{i_k}|\ \D_S(\bar\pi^{-1})\subseteq \bar B\}}
q^{{\ell}_S(\bar\pi)}t^{\bar\epsilon_S(\bar\pi)}\cdot
\left({n-1\brack i_k-1}_q+t_{i_k} q^{i}{n-1\brack i_k}_q\right).
$$

Now by a similar argument, this time applying Lemma~\ref{der7}(1), 

$$
\sum\limits_{\{\pi\in S_n|\ \D_S(\pi^{-1})\subseteq B\}} 
q^{\ell_S(\pi)}t^{\bar\ep_S(\pi)}=
$$
$$
\sum\limits_{\{\bar\pi\in S_{i_k}|\ \D_S(\bar\pi^{-1})\subseteq \bar B\}}
q^{{\ell}_S(\bar\pi)}t^{\bar\epsilon_S(\bar\pi)}\cdot
\left({n-1\brack i_k-1}_q+t_{i_k} q^{i}{n-1\brack i_k}_q\right),
$$
and the proof follows.
\qed

%

\bigskip

\noindent{\bf Proof of Theorem \ref{shuff.1}(2).}
\noindent
By the principle of inclusion and exclusion and Remark~\ref{rem5},
Theorem \ref{shuff.1}(2) is equivalent to the following statement :

For every subset $B\subseteq [n-1]$ 
$$
\sum\limits_{\{\sigma\in A_{n+1}|\ \D_A(\sigma^{-1})\subseteq B\}} 
q^{\rmaj_{A_{n+1}}(\sigma)}t^{\bar\ep_A(\sigma)}
=
$$
$$
\sum\limits_{\{\sigma\in A_{n+1}|\ \D_A(\sigma^{-1})\subseteq B\}} q^{\ell_A(\sigma)}t^{\bar\ep_A(\sigma)},
$$
By Proposition \ref{pro5} 
this part is reduced to Theorem \ref{shuff.1}(1).

\qed

\section{Appendix}
In this section we present another pair of statistics, leading to
a different analogue of MacMahon's Theorem.

\bigskip

\noindent 
For $1\le i< n$
define a map $h_i:S_n \longmapsto S_n$ as follows:
$$
h_i(\pi):= \cases
{s_i\pi, & \hbox{ if $i\in \D_S(\pi^{-1})$;}\cr
\pi, & \hbox{ if $i\not\in \D_S(\pi^{-1})$.}\cr}
$$
For every permutation $\pi\in S_n$ define
$$
\hat\ell_i (\pi):=\ell_S(h_i(\pi)),
$$
and
$$
\hat\maj_i(\pi):=\maj_S(h_i(\pi)).
$$
Then $\hat\ell_i$ and $\hat\maj_i$
are equi-distributed over the
even permutations in $S_n$ (i.e. over the alternating group $A_n$).

\begin{thm}\label{ap1}
Let $n\ge 2$, then
$$
\sum\limits_{\pi\in A_n}q^{\hat\ell_i(\pi)}=
\sum\limits_{\pi\in A_n}q^{\hat\maj_i(\pi)}=
\prod\limits_{i=3}^n (1+q+\dots+q^{i-1}).
$$
\end{thm}

\noindent{\bf Proof.}
By definition,
$$
Image(h_i)=\{\pi\in S_n\ |\ i\not\in Des_S(\pi^{-1})\}=
\{\pi\in S_n\ |\ \pi^{-1}\hbox{ is an\ $[n]\setminus \{i\}$-shuffle}\}.
$$
Also, for each $\sg\in Image(h_i)$,  
$h_i^{-1}(\sg)=\{\sg, s_i\sg\}$, and  exactly one element in the set $\{\sg, s_i\sg\}$ is even.

Thus, by Garsia-Gessel's Theorem (Theorem~\ref{p3}), 
$$
\sum\limits_{\pi\in A_n}q^{\hat\maj_i(\pi)}=
$$
$$
\sum\limits_{\{\pi\in S_n|\ \pi^{-1}\hbox{ is an\ $[n]\setminus \{i\}$-shuffle}\}}q^{\maj(\pi)}
={n\brack 2,1,\dots,1}_q=
\prod\limits_{i=3}^n (1+q+\dots+q^{i-1}),
$$
and similarly for $\hat\ell_i$.

\qed

\end{document}